\documentclass[11pt,openbib]{article}
\pdfoutput=1
\usepackage{amssymb}
\usepackage{amsmath}
\usepackage{amsfonts}
\usepackage[T1]{fontenc}
\usepackage{graphicx}

\newcommand{\rd}{\mathrm{d}}
\newcommand{\re}{\mathrm{e}}
\newcommand{\ri}{\mathrm{i}}
\newcommand{\R}{\mathbb{R}} 

\newcommand{\Z}{\mathbb{Z}}
\newcommand{\su}{\mathrm{su}}

\newcommand{\onehalf}{\mbox{$\frac{\scriptstyle 1}{\scriptstyle 2}$}} 
 
\newcommand{\lefthook}{\mbox{$\, \rule{8pt}{.5pt}\rule{.5pt}{6pt}\, \, $}}
\newcommand{\ttfrac}[2]{\mbox{$\frac{{\scriptstyle #1}}{{\scriptstyle #2}}$}}
\newcommand{\tttfrac}[2]{\mbox{$\frac{{\scriptscriptstyle #1}}{{\scriptscriptstyle #2}}$}}

\allowdisplaybreaks

\begin{document}

\title{Bohr-Sommerfeld-Heisenberg Quantization of the 2-dimensional Harmonic
Oscillator}
\author{Richard Cushman and J\k{e}drzej \'{S}niatycki}
\date{}

\maketitle

\begin{abstract}
We study geometric quantization of the harmonic oscillator in terms of a
singular real polarization given by fibres of the energy momentum map. \medskip  

\end{abstract}

\label{first}
\section{Introduction}

In the framework of geometric quantization we have formulated a quantization
theory based on the Bohr-Sommerfeld quantization rules \cite{bohr}, \cite{sommerfeld}. 
The aim of this theory is to provide an alternative approach
to quantization of a completely integrable system that incorporates the
Bohr-Sommerfeld spectrum of the defining set of commuting dynamical
variables of the system. The resulting theory closely resembles Heisenberg's
matrix theory and we refer to it as the Bohr-Sommerfeld-Heisenberg
quantization \cite{cushman-sniatycki12a}. \medskip

In this paper we give a detailed treatment of Bohr-Sommerfeld-Heisenberg
quantization of the $2$-dimensional harmonic oscillator. We also reduce the $S^1$ 
oscillator symmetry of the $2$-dimensional harmonic  
oscillator and a discrete ${\Z }_2$ symmetry. We quantize the $S^1$-reduced system and 
the ${\Z }_2 \times S^1$-reduced system. \medskip

We now give an overview of Bohr-Sommerfeld-Heisenberg quantization. 
Let $P$ be a smooth manifold with symplectic form $\omega $. Consider a
completely integrable system $(f_{1},...,f_{n},P,\omega )$, where $n$ is
equal to $\onehalf \dim P$, and 
$f_{1},...,f_{n}$ are Poisson commuting functions that are independent on 
an open dense subset $U$ of $P$. If all integral curves of the Hamiltonian vector fields of 
$f_{1},...,f_{n}$ are closed, then the joint level sets of $f_{1},...,f_{n}$
form a singular foliation of $(P,\omega )$ by tori. The restriction of this 
singular foliation to the open dense subset $U$
of $P$ is a regular foliation of $(U,\omega _{\mid U})$ by Lagrangian tori. \medskip

Each torus $T$ of the regular foliation has a neighbourhood $W$ in $P$ such that the
restriction of $\omega $ to $W$ is exact, see \cite[chapter IX]{cushman-bates}. 
In other words, $\omega _{\mid W}=\rd  \theta _{W}$. The modern version of Bohr-Sommerfeld quantization 
rule requires that for each generator $\Gamma _{i}$ of the fundamental
group of the $n$-torus $T,$ we have 
\begin{equation}
\int_{{\Gamma }_{i}}\theta _{W}=m_{i}h\mspace{10mu}\text{for}\mspace{5mu}%
i=1,...,n,  \label{Bohr-Sommerfeld}
\end{equation}%
where $m_{i}$ is an integer and $h$ is Planck's constant. It is easy to
verify that the Bohr-Sommerfeld conditions are independent of the choice of
the form $\theta _{W}$ satisfying $\rd  \theta _{W}=\omega _{\mid W}$. 
If $T$ is in the complement of $U$, then the Bohr-Sommerfeld conditions hold for 
$i=1, \ldots , \dim T$. Intrinsically, the Bohr-Sommerfeld
conditions are equivalent to the requirement that the connection on the 
prequantization line bundle $L$ when restricted to $T$ has
trivial holonomy group \cite{sniatycki}. \medskip

Let $S$ be the collection of all tori satisfying the Bohr-Sommerfeld
conditions. We refer to $S$ as the Bohr-Sommerfeld set of the integrable
system $(f_1, \ldots , f_n, P, \omega )$. Since the curvature form of $L$ is
symplectic, it follows that the complement of $S $ is open in $P$. Hence,
the representation space $\mathfrak{H}$ of geometric quantization of an
integrable system consists of distribution sections of $L$ supported on the
Bohr Sommerfeld set $S$. Since these distributional sections are covariantly
constant along the leaves of the foliation by tori, it follows that each torus $T\in S$ corresponds 
to a 1-dimensional subspace $\mathfrak{H} _{T}$
of $\mathfrak{H}$. We choose an inner product $(\, \, | \, \, )$ on $\mathfrak{H}$ so 
that the family $\{\mathfrak{H}_{T}\, | \, \mathrm{T} \in S \}$ consists of mutually
orthogonal subspaces. \medskip

In Bohr-Sommerfeld quantization, one assigns to each $n$-tuple of Poisson
commuting constants of motion $f=(f_{1},...,f_{n})$ on $P$ an $n$-tuple 
$({\mathbf{Q}}_{f_{1}},\ldots ,{\mathbf{Q}}_{f_{n}})$ of commuting quantum
operators ${\mathbf{Q}}_{f_{k}}$ for $1\leq k\leq n$ such that for each $n$%
-torus $T\in S$, the corresponding $1$-dimensional space ${\mathfrak{H}}_{T}$
of the representation space $(\mathfrak{H},(\,\,|\,\,))$ is an eigenspace
for each ${\mathbf{Q}}_{f_{k}}$ for $1\leq k\leq n$ with eigenvalue ${f_k}|_T$. 
For any smooth function $F\in C^{\infty }({\R }^{n}),$
the composition $F(f_{1},...,f_{n})$ is quantizable. The operator $%
{\mathbf{Q}}_{F(f_{1},...,f_{n})}$ acts on each $\mathfrak{H}_{T}$ by
multiplication by $F(f_{1},...,f_{n})_{\mid T}$. \medskip

We assume that there exist global action-angle variables $(A_{i},\varphi _{i})$ on $U$ such that 
$\omega _{\mid U}= \rd  (\sum_{i=1}^{n}A_{i} \rd  \varphi _{i})$. Therefore, we can replace 
$\theta _{W}$ by the $1$-form $\sum^n_{i=1} A_i \, \rd  {\varphi }_i$
in equation (\ref{Bohr-Sommerfeld}) and rewrite the Bohr sommerfeld
conditions as 
\begin{equation}
\int_{{\Gamma }_{i}}A_{i}\, \rd  {\varphi }_{i}=m_{i}h\mspace{10mu}\text{for each}\mspace{5mu}i=1,...,n.
\label{Bohr-Sommerfeld2}
\end{equation}%
Since the actions $A_{i}$ are independent of the angle variables, we can
perform the integration and obtain 
\begin{equation}
A_{i}=m_{i}\hbar \mspace{10mu}\text{for each}\mspace{5mu}i=1,...,n.
\label{Bohr-Sommerfeld1}
\end{equation}%
where $\hbar =h/2\pi $. For each multi-index $\mathbf{m}=(m_{1},...,m_{n})$,
we denote by $T_{\mathbf{m}}$ the torus satisfying the Bohr-Sommerfeld
conditions with integers $m_{1},...,m_{n}$ on the right hand side of
equation (\ref{Bohr-Sommerfeld1}). \medskip

Let $S_{U}$ be the restriction of the Bohr-Sommerfeld set $S$ to the open
neighbourhood $U$ of $P$ on which the functions $f_{1},...,f_{n}$ are
independent. Consider a subspace ${\mathfrak{H}}_{U}$ of $\mathfrak{H}$
given by the direct sum of $1$-dimensional subspaces ${\mathfrak{H}}_{T_{\mathbf{m}}}$ 
corresponding to Bohr-Sommerfeld tori $T_{\mathbf{m}}\in S_{U}$. Let ${\mathbf{e}}_{\mathbf{m}}$ be a basis vector of $\mathfrak{H}_{T_{\mathbf{m}}}$. Each ${\mathbf{e}}_{\mathbf{m}}$ is a joint eigenvector of
the commuting operators $({\mathbf{Q}}_{A_{1}},\ldots ,{\mathbf{Q}}_{A_{n}})$ corresponding to the eigenvalue 
$(m_{1}\hbar ,\ldots ,m_{n}\hbar) $. The vectors $({\mathbf{e}}_{\mathbf{m}})$ form an orthogonal basis in 
${\mathfrak{H}}_{U}$. Thus, 
\begin{equation}
\left( {\mathbf{e}}_{\mathbf{m}}\mid {\mathbf{e}}_{{\mathbf{m}}^{\prime}}\right) =0\mspace{10mu}\text{if}\mspace{5mu}\mathbf{m}\neq {\mathbf{m}}^{\prime }.  
\label{product1}
\end{equation}

For each $i=1,...,n$, introduce an operator ${\mathbf{a}}_{i}$ on ${\mathfrak{H}}_{U}$ such that 
\begin{equation}
{\mathbf{a}}_{i}{\mathbf{e}}_{(m_{1},..,,m_{i-1},m_{i},m_{i+1},...,m_{n})}=
{\mathbf{e}}_{(m_{1},...m_{i-1},m_{i}-1,m_{i+1},...,m_{n})\text{.}}  
\label{5one}
\end{equation}%
In other words, the operator ${\mathbf{a}}_{i}$ shifts the joint
eigenspace of $({\mathbf{Q}}_{A_{1}},\ldots ,{\mathbf{Q}}_{A_{n}})$
corresponding to the eigenvalue $(m_{1}\hbar ,\ldots ,m_{n}\hbar )$ to the
joint eigenspace of $({\mathbf{Q}}_{A_{1}},\ldots ,$ ${\mathbf{Q}}_{A_{n}})$
corresponding to the eigenvalue $(m_{1}\hbar ,$ $\ldots ,m_{i-1}\hbar $, $%
(m_{i}-1)\hbar $, $m_{i+1}\hbar ,\ldots ,$ $m_{n}\hbar )$. Let $%
{\mathbf{a}}_{i}^{\dagger }$ be the adjoint of ${\mathbf{a}}_{i}$.
Equations (\ref{product1}) and (\ref{5one}) yield 
\begin{equation}
{\mathbf{a}}_{i}^{\dagger }{\mathbf{e}}_{(m_{1},..,,m_{i-1},m_{i},m_{i+1},...,m_{n})} =
{\mathbf{e}}_{(m_{1},...m_{i-1},m_{i}+1,m_{i+1},...,m_{n})}  
\label{6}
\end{equation}%
We refer to the operators ${\mathbf{a}}_{i}$ and ${\mathbf{a}}_{i}^{\dagger }$ as \emph{shifting} operators.\footnote{In representation theory, shifting operators are called ladder operators.
The corresponding operators in quantum field theory are called the creation
and annihilation operators.} For every $i=1,...,n$, we have 
\begin{equation}
\lbrack {\mathbf{a}}_{i}, {\mathbf{Q}}_{A_{j}}]=\hbar \,{\mathbf{a}}_{i} \, \delta _{ij},  
\label{commutation}
\end{equation}%
where $\delta _{ij}$ is Kronecker's \ symbol that is equal to $1$ if $i=j$
and vanishes if $i\neq j$. Taking the adjoint, of the preceding equation, we
get 
\begin{equation*}
\lbrack {\mathbf{a}}_{i}^{\dagger }, {\mathbf{Q}}_{A_{j}}]=-\hbar \, {\mathbf{a}}_{i}^{\dagger }\, \delta _{ij}.
\end{equation*}

The operators ${\mathbf{a}}_{i}$ and ${\mathbf{a}}_{i}^{\dagger }$ are
well defined in the Hilbert space $\mathfrak{H}$. In \cite%
{cushman-sniatycki12a}, we have interpreted them as the result of
quantization of apropriate functions on $P$. We did this as follows. We look
for smooth complex-valued functions $h_{i}$ on $P$ satisfying the Poisson
bracket relations 
\begin{equation}
\{h_{j},A_{i}\}= - \ri \, \delta _{ij}h_{j}.  
\label{bracket}
\end{equation}%
The Dirac quantization condition  
\begin{equation}
\lbrack {\mathbf{Q}}_{f},{\mathbf{Q}}_{h}]= \ri \hbar \, {\mathbf{Q}}_{\{f,h\}}  
\label{Dirac}
\end{equation}%
implies that we may interpret the operator ${\mathbf{a}}_{j}$ as the
quantum operator corresponding to $h_{j}$. In other words, we set 
${\mathbf{a}}_{j}={\mathbf{Q}}_{h_{j}}$. This choice is consistent with
equation (\ref{bracket}) because (\ref{commutation}) yields%
\begin{equation*}
\lbrack {\mathbf{Q}}_{h_{j}},{\mathbf{Q}}_{A_{i}}]= \ri \hbar \, {\mathbf{Q}}_{\{h_{j},A_{i}\}}= 
\ri \hbar \,{\mathbf{Q}}_{-\ri \delta _{ij}h_{j}}= 
{\delta}_{ij}\hbar \,{\mathbf{Q}}_{h_{j}}\text{.}
\end{equation*}

Since $\omega _{\mid U}=\sum_{i=1}^{n}\, \rd  A_{i}\wedge \rd  \varphi _{i}$, it
follows that the Poisson bracket of ${\re }^{-\ri \varphi _{j}}$ and $A_{i}$ is 
\begin{equation}
\{ {\re}^{-\ri \varphi _{j}},A_{i}\}=X_{A_{i}}{\re}^{-\ri \varphi _{j}}=%
\tfrac{\partial }{\partial \varphi _{i}}{\re}^{-\ri \varphi _{j}}= 
-\ri \, {\delta }_{ij}{\re}^{-\ri \varphi _{j}}.  
\label{bracket2}
\end{equation}%
Comparing equations (\ref{bracket}) and (\ref{bracket2}) we see that we may
make the following identification ${\mathbf{a}}_{i}=
{\mathbf{Q}}_{{\re}^{-\ri \varphi _{i}}}$ and ${\mathbf{a}}_{i}^{\dagger }=
{\mathbf{Q}}_{{\re }^{\ri \varphi _{i}}}$ for $i=1,...,n$. Clearly, the functions 
$h_{j}={\re }^{\pm \ri \varphi _{j}}$ are not uniquely defined by equation (\ref{bracket}). We 
can multiply them by arbitrary functions that commute 
with all actions $A_{1},\ldots ,A_{n}$. Hence, there is a choice
involved. We shall use this freedom of choice to satisfy consistency
requirements on the boundary of $U$ in $P$. \medskip

\section{The classical theory}

We describe the geometry of the $2$-dimensional harmonic
oscillator. \medskip

The configuration space of the $2$-dimensional harmonic oscillator is ${\R}^2$ with coordinates 
$x=(x_1,x_2)$. The phase space is $T^\ast {\R }^2 
= {\R }^4$ with coordinates $(x,y)=(x_1,x_2,y_1,y_2)$.
On $T^{\ast }{\R  }^2$ the canonical $1$-form is ${\Theta }^{\scriptscriptstyle \vee }= y_1 \, \rd  x_1 + y_2 \, \rd  x_2
= \langle y, \rd  x \rangle $. The symplectic form on $T^{\ast}{\R  }^2$ is the closed nondegenerate 
$2$-form ${\omega }^{{\scriptscriptstyle \vee }}= \rd  {\Theta }^{\scriptscriptstyle \vee } = \rd  y_1 \wedge \rd  x_1+ \rd  y_2 \wedge \rd  x_2$ whose matrix representation is 
\begin{equation}
\omega = {\begin{pmatrix}
\rd  x \\ 
\rd  y
\end{pmatrix}
}^T 
\begin{pmatrix}
0 & -I_2 \\ 
I_2 & 0
\end{pmatrix}
\begin{pmatrix}
\rd  x \\ 
\rd  y
\end{pmatrix} .
\label{eq-sec1dot}
\end{equation}
Given a smooth function $f \in C^{\infty}(T^{\ast }{\R  }^2)$, its Hamiltonian vector field is 
$({\omega }^{\flat })^{-1}\rd f = \langle \frac{\partial f}{\partial y} , \frac{\partial }{\partial x} \rangle - 
\langle \frac{\partial f}{\partial y} , \frac{\partial }{\partial x} \rangle $. Let $(x,y) \in T^{\ast }{\R  }^2$. For 
every $f$, $g \in C^{\infty}(T^{\ast }{\R  }^2)$ their Poisson bracket is $\{ f, g \} (x,y) = 
\omega (x,y)\big( X_g(x,y), X_f(x,y) \big) $, which is a smooth function on $T^{\ast }{\R  }^2$. 
The structure matrix {\tiny $\left( 
\begin{array}{c|c}
\{ x_i, x_j \} & \{ x_i, y_j \} \\ \hline
\{ y_i, x_j \} & \{ y_i, y_j \}
\end{array}
\right) $} of the Poisson bracket $\{ \, \, , \, \, \} $ on $C^{\infty}(T^{\ast }{\R  }^2)$ is 
{\tiny $\begin{pmatrix}
0 & I_2 \\ 
-I_2 & 0
\end{pmatrix}
$}. \medskip

The Hamiltonian function of the $2$-dimensional harmonic oscillator is 
\begin{equation}
E^{{\scriptscriptstyle \vee }}:T^{\ast }{\R }^{2}\rightarrow \R :(x,y)\mapsto 
\onehalf (x_{1}^{2}+y_{1}^{2})+
\onehalf (x_{2}^{2}+y_{2}^{2}).
\label{eq-sec1dagger}
\end{equation}%
The corresponding Hamiltonian vector field $X_{E^{{\scriptscriptstyle \vee }}}=\langle X_{1},\frac{\partial }{\partial x}\rangle +\langle X_{2},\frac{\partial }{\partial y}\rangle $ can be computed using $-\rd  E^{{\scriptscriptstyle \vee }}=X_{E^{{\scriptscriptstyle \vee }}} \lefthook \, \omega =\langle X_{2},\rd  x\rangle -\langle X_{1},\rd  y\rangle $. We get 
\begin{equation}
X_{1}=\mbox{$\frac{{\scriptstyle \partial E^{{\scriptscriptstyle \vee }}}}{{\scriptstyle \partial y}}$}%
\quad \text{and}\quad X_{2}=-\mbox{$\frac{{\scriptstyle \partial E^{{\scriptscriptstyle \vee }}}}{{\scriptstyle \partial x}}$}.
\end{equation}%
Therefore the equations of motion of the harmonic oscillator are 
\begin{equation}
\dot{x}=y\quad \text{and}\quad \dot{y}=-x.
\label{eq-sec1dot}
\end{equation}%
The solution to the above equations is the one parameter family of
transformations 
\begin{equation}
\phi _{t}^{E^{{\scriptscriptstyle \vee }}}(x,y)=A(t)%
\begin{pmatrix}
x \\ 
y%
\end{pmatrix}%
\,=\,%
\begin{pmatrix}
(\cos t)\,I_{2} & -(\sin t)\,I_{2} \\ 
(\sin t)\,I_{2} & (\cos t)\,I_{2}%
\end{pmatrix}%
\begin{pmatrix}
x \\ 
y%
\end{pmatrix}%
\end{equation}%
This defines an $S^{1}$ action, called the \emph{oscillator symmetry},
on $T^{\ast }{\R }^{2}$, which is a map from $\R $ to 
$\mathrm{Sp} (4,\R )$ that sends $t$ to the $4\times 4
$ symplectic matrix $A(t)$, which is periodic of period $2\pi $. \medskip 

Since $E^{{\scriptscriptstyle \vee }}$ is constant along the integral curves of $X_E$, the manifold 
\begin{equation}
(E^{{\scriptscriptstyle \vee }})^{-1}(e) = \{(x,y)\in {\R }^4 \, \mathop{\rule[-4pt]{.5pt}{13pt}\, }%
\nolimits x^2+y^2=2e, e>0 \} = S_{\sqrt{2e}}^3,
\end{equation}
which is a $3$-sphere of radius $\sqrt{2e}$, is invariant under the flow of $X_{E^{{\scriptscriptstyle \vee }}}$. \medskip

The configuration space ${\R }^{2}$ is invariant under the $S^{1}$
action $S^{1}\times {\R }^{2}\rightarrow {\R }^{2}:(t,x)\mapsto R_{t}x$, 
where $R_{t}$ is the matrix {\tiny $%
\begin{pmatrix}
\cos t & -\sin t \\ 
\sin t & \hfill \cos t
\end{pmatrix}
$}. This map lifts to a symplectic action $\Psi _{t}$ of $S^{1}$ on phase
space $T^{\ast }{\R }^{2}$ that sends $(x,y)$ to ${\Psi }_{t}(x,y)=(R_{t}x,R_{t}y)$. 
The infinitesimal generator of this action is 
\begin{equation}
Y(x,y)=\mbox{${\displaystyle \frac{\rd  }{\rd  t}}
\rule[-10pt]{.5pt}{25pt} \raisebox{-10pt}{$\, {\scriptstyle t=0}$}$}\hspace{%
-8pt}\Phi _{t}(x,y)=(-x_{2},x_{1},-y_{2},y_{1}).
\end{equation}%
The vector field $Y$ is Hamiltonian corresponding to the Hamiltonian
function 
\begin{equation}
L^{{\scriptscriptstyle \vee }}(x,y)=\langle (y_1,y_2),(x_{2},-x_{1})\rangle =x_{1}y_{2}-x_{2}y_{1},
\label{eq-sec1neweq}
\end{equation}%
that is, $Y=X_{L^{{\scriptscriptstyle \vee }}}$. $L^{{\scriptscriptstyle \vee } }$ is readily recognized as the angular momentum. The
integral curve of the vector field $X_{L^{{\scriptscriptstyle \vee }}}$ on ${\R }^{4}$ starting at 
$(x,y)=(x_{1},x_{2},y_{1},y_{2})$ is 
\begin{equation*}
s\mapsto {\varphi }_{s}^{L^{{\scriptscriptstyle \vee }}}(x,y)=%
\mbox{\footnotesize $\left( \begin{array}{rcrc}
\cos s & -\sin s & 0 & 0  \\
\sin s & \cos s & 0 & 0 \\
0 & 0 &  \cos s & -\sin s \\
0 & 0 &\sin s & \cos s \end{array}  \right) $}\,%
\mbox{\footnotesize 
$\begin{pmatrix}x \\ y \end{pmatrix} $},
\end{equation*}%
which is periodic of period $2\pi $. The Hamiltonian of the harmonic
oscillator is an integral of $X_{L^{{\scriptscriptstyle \vee } }}$. The $S^{1}$ symmetry $\Psi $ of the
harmonic oscillator is called the \textit{angular momentum} symmetry. \medskip 

Consider the completely integrable system $(E^{{\scriptscriptstyle \vee } },L^{{\scriptscriptstyle \vee } },{\R }^{4},{\omega }^{{\scriptscriptstyle \vee } })$
where $E^{{\scriptscriptstyle \vee } }$ (\ref{eq-sec1dagger}) and $L^{{\scriptscriptstyle \vee } }$ (\ref{eq-sec1neweq}) are the energy
and angular momentum of the $2$-dimensional harmonic oscillator. The
Hamiltonian vector fields 
\begin{equation*}
X_{E^{{\scriptscriptstyle \vee }}}=%
\mbox{\footnotesize $y_1\frac{\partial }{\partial x_1} + y_2\frac{\partial }{\partial x_2} - 
x_1\frac{\partial }{\partial y_1} - x_2 \frac{\partial }{\partial y_2}$}%
\mspace{10mu}\text{and}\mspace{10mu}X_{L^{\scriptscriptstyle \vee }}=%
\mbox{\footnotesize
$-x_2\frac{\partial }{\partial x_1} + x_1\frac{\partial }{\partial x_2} -
y_2\frac{\partial }{\partial y_1}  + y_1 \frac{\partial }{\partial y_2}$}
\end{equation*}%
corresponding to the Hamiltonians $E^{{\scriptscriptstyle \vee } }$ and $L^{{\scriptscriptstyle \vee } }$, respectively, define the
generalized distribution 
\begin{equation}
D:{\R }^{4}\rightarrow T{\R }^{4}:(x,y)\mapsto D_{(x,y)}=%
{\rm span} \{X_{E^{{\scriptscriptstyle \vee }}}(x,y),\,X_{L^{{\scriptscriptstyle \vee }}}(x,y)\}.  
\label{eq-polone}
\end{equation}

Before we can investigate the geometry of $D$ (\ref{eq-polone}) we will need
some very detailed geometric information about the $2$-torus action 
\begin{equation}
\Phi : T^2 \times {\R  }^4 \rightarrow {\R  }^4:
\big( (t_1,t_2), (x,y) \big) \mapsto {\varphi }^{E^{{\scriptscriptstyle \vee }}}_{t_1} \raisebox{0pt}{$\scriptstyle%
\circ \, $} {\varphi }^{L^{{\scriptscriptstyle \vee } }}_{t_2}(x,y)  \label{eq-poltwo}
\end{equation}
generated by the $2\pi $ periodic flows ${\varphi }^{E^{{\scriptscriptstyle \vee } }}_{t_1}$ and ${\varphi }^{L^{{\scriptscriptstyle \vee }}}_{t_2}$ of the vector fields $X_{E^{{\scriptscriptstyle \vee } }}$ and $X_{L^{{\scriptscriptstyle \vee } }}$, respectively. Using
invariant theory we find the space $V = {\R  }^4/T^2$ of orbits of
the $T^2$-action $\Phi $ (\ref{eq-poltwo}) as follows. The algebra of $T^2$%
-invariant polynomials on ${\R  }^4$ is generated by ${\sigma }_1 = %
\onehalf (y^2_1+x^2_1 +y^2_2+x^2_2)$
and ${\sigma }_2 = x_1y_2-x_2y_1$. Let ${\sigma }_3 = \onehalf (y^2_1+x^2_1 -y^2_2-x^2_2)$ and ${\sigma 
}_4 = x_1y_1 +x_2y_2$. Then ${\sigma }^2_1 = 
{\sigma }^2_2+{\sigma }^2_3 +{\sigma }^2_4 \ge {\sigma }^2_2$. From the fact 
that ${\sigma }_1\ge 0$, we
get ${\sigma }_1 \ge |{\sigma }_2|$. So the $T^2$ orbit space $V$ is the
semialgebraic variety 
\begin{displaymath}
\{ ( {\sigma }_1 , {\sigma }_2) \in {\R  }^2 \, %
\mathop{\rule[-4pt]{.5pt}{13pt}\, }\nolimits \, |{\sigma }_2 | \le {\sigma }_1, \, \, {\sigma }_1 \ge 0 \}. 
\end{displaymath}
Moreover $V$ is a differential space with its differential structure given by 
the space ${C^{\infty}({\R  }^4)}^{T^2}$ of smooth $T^2$-invariant functions on ${\R  }^4$. That is the 
differential structure of $V$ is the space of smooth functions in ${\sigma }_1$ and ${\sigma }_2$. Note that the $T^2$-orbit map $\pi : {\R  }^4
\rightarrow V \subseteq {\R  }^2:(x,y) \mapsto \big( {\sigma }_1(x,y), {\sigma }_2 (x,y) \big)$ 
is just the energy momentum mapping 
\begin{equation*}
(\mathcal{EM})^{{\scriptscriptstyle \vee } }:{\R  }^4 \rightarrow {\R  }^2:
(x,y) \mapsto \big( E^{{\scriptscriptstyle \vee } }(x,y), L^{{\scriptscriptstyle \vee } } (x,y) \big)
\end{equation*}
of the $2$-dimensional harmonic oscillator. Abstractly, the Whitney
stratification of $V$ is ${\coprod}_{0 \le i \le 2 }V_i $ where  
\begin{equation*}
V_2 = V  \setminus V_{\mathrm{sing}}, \, V_1 = V_{\mathrm{sing}} \setminus (V_{\mathrm{sing}})_{\mathrm{sing}}, 
\, \, \mathrm{and} \, \, V_0 = (V_{\mathrm{sing}})_{\mathrm{sing}} .
\end{equation*}
Here $W_{\mathrm{sing}}$ is the semialgebraic variety of singular points of
the semialgebraic variety $W$. We have $V_2 = \{ ({\sigma }_1, {\sigma }_2 ) \in {%
\R  }^2 \,  | \, |{\sigma }_2 | < {\sigma }_1 \}$, $V_1 =\{ ({\sigma }_1, {\sigma }_2 ) \in {\R  }^2 \, | 
\, 0 < |{\sigma }_2 | = {\sigma }_1 \} $, and 
$V_0 = \{ ({\sigma }_1, {\sigma }_2 ) \in {\R  }^2 \, | \, {\sigma }_1 = {\sigma }_2=0 \}$.  \medskip

For $j=0,1,2$ let $U_j = {\mathcal{EM}}^{-1}(V_j)$. It is straightforward to
see that the isotropy group $T^2_{(x,y)}$ at $(x,y)$ of the $T^2$-action $\Phi $ (\ref{eq-poltwo}) is 
\begin{equation*}
T^2_{(x,y)} = \left\{ 
\begin{array}{cl}
\{ e \} , & \mbox{if $(x,y) \in U_2$} \\ 
S^1 & \mbox{if $(x,y) \in U_1 $} \\ 
T^2, & \mbox{if $(x,y) \in U_0$.}%
\end{array}
\right.
\end{equation*}
Thus $\{ U_j, \, j=0,1,2 \} $ is the stratification of ${\R  }^4$ by
orbit type, which is equal to $T^2$-symmetry type where all the isotropy
groups are equal, since $T^2$ is abelian. Because $V_j$ is the image of $U_j$
under the energy momentum map, we see that $U_j$ is the subset of ${\R  }^4$ where the rank of 
$D\mathcal{EM}$ is $j$. Since 
\begin{equation*}
{\omega }^{\flat}(x,y):T_{(x,y)}{\R  }^4 \rightarrow 
T^{\ast }_{(x,y)}{\R  }^4:v_{(x,y)} \mapsto \{ w_{(x,y)}  \mapsto \omega (x,y)(v_{(x,y)}, w_{(x,y)}) \} 
\end{equation*}
is bijective, it follows that 
\begin{align}
\dim D_{(x,y)} & = \dim \mathrm{span} \{ \rd  E(x,y), \, \rd   L(x,y) \} = 
\mathrm{rank}\, D\mathcal{EM}(x,y).  
\label{eq-polthree}
\end{align}
Therefore $U_j$ is the union of $j$-tori for $j=0,1,2$ each of which is an
energy momentum level set of $(\mathcal{EM})^{\scriptscriptstyle \vee }$. \medskip

We turn to investigating the generalized distribution $D$ (\ref{eq-polone}).
Since $\omega (X_{E^{{\scriptscriptstyle \vee } }},$ $X_{L^{{\scriptscriptstyle \vee }}}) = 
\{ L^{{\scriptscriptstyle \vee } }, E^{{\scriptscriptstyle \vee } } \} $, we see that for every $(x,y) \in 
{\R  }^4$ we have ${\omega }(x,y)D_{(x,y)} = 0$. Thus every subspace 
$D_{(x,y)}$ in the distribution $D$ is an $\omega $-isotropic subspace of the
symplectic vector space $\big( T_{(x,y)}{\R  }^4, \omega (x,y) \big)$. 
If $(x,y) \in U_j$ then $\dim D_{(x,y)} = j$. So if $(x,y) \in U_2$, then $D_{(x,y)}$ is 
a Lagrangian subspace of $\big( T_{(x,y)}{\R  }^4,
\omega (x,y) \big)$. If $u_{(x,y)}$, $v_{(x,y)} \in D_{(x,y)}$, then
thinking of $u_{(x,y)}$ and $v_{(x,y)}$ as vector fields on $T{\R  }^4 $ 
we obtain $[u_{(x,y)}, v_{(x,y)}] =0$, because their flows commute.
Thus $D $ is an involutive generalized distribution, which is called the
\emph{energy-momentum polarization} of the symplectic manifold $({\R  }^4,
\omega )$. The leaf (= integral manifold) $L_{(x,y)}$ of the distribution $D$
through the point $(x,y) \in {\R  }^4$ is diffeomorphic to $T^2/T^2_{(x,y)}$. In particular we have 
\begin{equation*}
L_{(x,y)} = \left\{ 
\begin{array}{cl}
T^2, & \mbox{if $(x,y) \in U_2$} \\ 
T^1 = S^1, & \mbox{if $(x,y) \in U_1$} \\ 
T^0= \mathrm{pt}, & \mbox{if $(x,y) \in U_0$.}
\end{array}
\right.
\end{equation*}
Thus $D$ is a generalized toric distribution. Note that the leaves of $D$,
each of which is an energy momentum level set of $(\mathcal{EM})^{{\scriptscriptstyle \vee } }$, foliate the
strata $U_j$ for $j=0,1,2$. Thus the space of leaves ${\R  }^4/D$ of $D$ is 
the $T^2$ orbit space $V$. \medskip 

\noindent \hspace{1in}\begin{tabular}{l}
\setlength{\unitlength}{.5pt} 
\includegraphics[width=200pt]{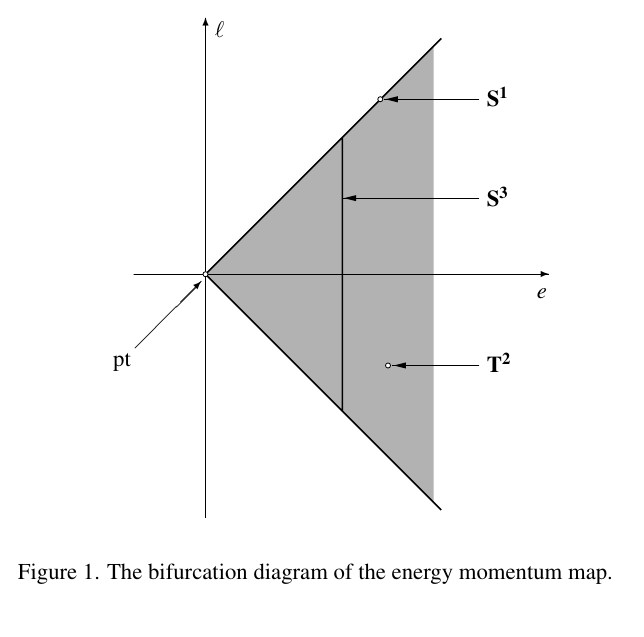}
\end{tabular}%
\medskip

We collect all the geometric information about the leaves of the generalized
distribution $D$ in the bifurcation diagram of the energy momentum map of
the $2$-dimensional harmonic oscillator above. The range of this map is
indicated by the shaded region. \medskip

We now find action-angle coordinates for the $2$-dimensional harmonic
oscillator. Consider the $T^2$ action on $(T^{\ast }{\R }^2, {\omega }^{\scriptscriptstyle \vee })$ generated by the flows of the Hamiltonian vector fields $X_{E^{{\scriptscriptstyle \vee } }}$ of energy $E^{{\scriptscriptstyle \vee } }$ (\ref{eq-sec1dagger}) and $X_{L^{{\scriptscriptstyle \vee } }}$ of angular momentum $L^{{\scriptscriptstyle \vee } }$ (\ref{eq-sec1neweq}), specifically, 
\begin{displaymath}
\Phi : T^2 \times T^{\ast }{\R }^2 \rightarrow T^{\ast }{\R }^2:\big( (t,s), (x,y) \big) \mapsto 
{\varphi }^{E^{{\scriptscriptstyle \vee } }}_t {\scriptstyle \circ} {\varphi }^{L^{{\scriptscriptstyle \vee } }}_s(x,y).
\end{displaymath}
Consider the change coordinates by 
\begin{equation} 
\psi :T^{\ast }{\R }^2 \rightarrow T^{\ast }{\R }^2: \mbox{\small $\begin{pmatrix} 
x \\ y \end{pmatrix}$} \longmapsto \mbox{\footnotesize $\begin{pmatrix} 
\xi \\ \eta  \end{pmatrix}$} =\mbox{\footnotesize $\mbox{$\frac{\scriptstyle 1}{\scriptstyle \sqrt{2}}$}\, 
\begin{pmatrix}  1 & 0 & 0 & 1 \\ -1 & 0 & 0 & 1 \\ 0 & -1 & 1 & 0 \\ 0 & -1 & -1 & 0 \end{pmatrix} $}\, \mbox{\small $\begin{pmatrix} x \\ y \end{pmatrix}$.}
\label{eq-sec2newone}
\end{equation}
Then $\psi $ is symplectic, that is $\omega  = {\psi }^{\ast } {\omega }^{{\scriptscriptstyle \vee } }= \rd {\eta }_1 \wedge \rd {\xi }_1 + 
\rd {\eta }_2 \wedge \rd {\xi }_2$. Moreover, it diagonalizes the Hamiltonians for $E^{{\scriptscriptstyle \vee } }$ and $L^{{\scriptscriptstyle \vee } }$, that is, 
in the new coordinates $E^{{\scriptscriptstyle \vee } }$ becomes 
\begin{equation}
E(\xi , \eta ) = {\psi }_{\ast }E^{{\scriptscriptstyle \vee } } (\xi , \eta ) = \onehalf  \, ( {\xi }^2_1 +{\eta }^2_1 +{\xi }^2_2 +{\eta }^2_2) ;
\label{eq-sec2newzeroa}
\end{equation}
while $L^{{\scriptscriptstyle \vee } }$ becomes 
\begin{equation}
L(\xi , \eta ) = {\psi }_{\ast }L^{{\scriptscriptstyle \vee }} (\xi , \eta ) = \onehalf  \, ( {\xi }^2_1 +{\eta }^2_1 - {\xi }^2_2 -{\eta }^2_2) .
\label{eq-sec2newzerob}
\end{equation}
Let 
\begin{subequations}
\begin{align}
A_1 & = \onehalf  \, ({\xi }^2_1 +{\eta }^2_1) = 
\onehalf  \big( E(\xi , \eta ) + L(\xi , \eta ) \big) \ge 0 
\label{eq-sec2newonea} \\
A_2 & = \onehalf  \, ({\xi }^2_2 +{\eta }^2_2) = 
\onehalf  \big( E(\xi , \eta ) - L(\xi , \eta ) \big) \ge 0 
\label{eq-sec2newoneb} 
\end{align}
\end{subequations}
and let 
\begin{equation}
A: T^{\ast }{\R }^2 \rightarrow ({\R }_{\ge 0})^2 \subseteq {\R }^2: (\xi , \eta ) \mapsto \big( A_1(\xi , \eta ), A_2(\xi ,\eta ) \big) .
\label{eq-sec2newtwo}
\end{equation}
Then the flow of $X_{A_1}$ and $X_{A_2}$ defines a $T^2$-action on $A^{-1}(({\R }_{> 0})^2)$ 
given by 
\begin{displaymath}
\begin{array}{l}
T^2 \times A^{-1}(({\R }_{> 0})^2) \subseteq T^{\ast }{\R }^2 \rightarrow 
A^{-1}(({\R }_{> 0})^2)\subseteq T^{\ast }{\R }^2 : 
\big( (t, s), (\xi , \eta ) \big) \longmapsto \\
\rule{0pt}{30pt} \hspace{.25in} {\varphi }^{E} {\scriptstyle \circ } {\varphi }^{L} (\xi ,\eta ) = 
\mbox{\footnotesize $\begin{pmatrix} \cos t & 0 & \sin t & 0 \\ 0 & 0 & 0& 0 \\
-\sin t & 0 & \cos t & 0 \\ 0 & 0 & 0 & 0 \end{pmatrix}$} \, 
\mbox{\footnotesize $\begin{pmatrix} 0 & 0 & 0 & 0 \\
0 & \cos s & 0 & \sin s  \\ 0 & 0 & 0& 0 \\
0 & -\sin s & 0 & \cos s \end{pmatrix}$} \, \mbox{\footnotesize $\begin{pmatrix} \xi \\ \eta \end{pmatrix}$, }
\end{array}
\end{displaymath}
which is periodic of period $2\pi $ on every $2$-torus $T^2_{a_1, a_2} = A^{-1}_1(a_1) \cap A^{-1}_2(a_2)$, 
where $(a_1, a_2) \in ({\R }_{> 0})^2$. So the functions $A_1$ (\ref{eq-sec2newonea}) and $A_2$ 
(\ref{eq-sec2newoneb}) are actions and the map $A$ (\ref{eq-sec2newtwo}) is the action map. The angle 
corresponding to the action $A_1$ is ${\vartheta }_1 = t = {\tan }^{-1}$ 
$\ttfrac{{\eta }_1}{{\xi }_1} $; while the angle corresponding to the action $A_2$ is ${\vartheta }_2 = s = 
{\tan }^{-1}$$\ttfrac{{\eta }_2}{{\xi }_2} $. A  calculation shows that 
$\omega = \rd A_1 \wedge \rd {\vartheta }_1 + \rd A_2 \wedge \rd {\vartheta }_2$. Thus $(A_1, A_2, 
{\vartheta }_1, {\vartheta}_2)$ are action-angle coordinates on $\big( A^{-1}(({\R }_{> 0})^2), 
{\omega }|_{A^{-1}(({\R }_{> 0})^2)} \big) $. When $A_2 =0$ on 
$\partial {A^{-1}(({\R }_{\ge 0})^2)}$ action-angle coordinates on 
$A^{-1}({\R }_{> 0} \times \{ 0 \}) $ are $(A_1, 0, {\vartheta }_1, 0)$; while when $A_1 =0$ on 
$\partial {A^{-1}(({\R }_{\ge 0})^2)}$ action-angle coordinates on 
$A^{-1}(\{ 0 \} \times {\R }_{> 0}) $ are $(0,A_2, 0, {\vartheta }_2)$. 

\section{Bohr-Sommerfeld-Heisenberg quantization}

In this section we carry out the Bohr-Sommerfeld-Heisenberg quantization of the $2$-dimensional 
harmonic oscillator. \medskip 

First we look at its Bohr-Sommerfeld quantization. This amounts to quantizing the action functions, namely, 
setting 
\begin{equation}
A_1 = m\hbar = a_1 \, \, \, \mathrm{and} \, \, \, A_2 = n\hbar = a_2 , 
\label{eq-secn3newone}
\end{equation}
where $(m,n) \in ({\Z }_{\ge 0})^2$. Next associate a vector ${\mathbf{e}}_{m,n}$ to each 
Bohr-Sommerfeld torus $A^{-1}(a_1, a_2)$ where 
$(a_1,a_2)$ lie in $({\R }_{\ge 0})^2$, which is the image of the action mapping $A$. In other words, 
$(a_1,a_2) =  \big( \onehalf \hbar (e+\ell ), 
\onehalf  \hbar (e-\ell ) \big) $, where  $\hbar (e, \ell ) = \hbar (m+n, m-n)$ lies in the image $V$ of the energy momentum mapping $(\mathcal{EM})^{\scriptscriptstyle \vee }$ of the harmonic oscillator. Let $\mathfrak{H}=\mathrm{span}\, \{ {\mathbf{e}_{m,n} \, |  \, (m,n} \in ({\Z }_{\ge 0})^2 \} $. 
Because the Bohr-Sommerfeld tori are disjoint, we may choose a hermitian inner product 
$\langle \, \, , \, \, \rangle $ on $\mathfrak{H}$ so that $\langle {\mathbf{e}}_{m,n}, {\mathbf{e}}_{m',n'} \rangle  = 
{\delta }_{(m,n), (m',n')} = {\delta }_{m,m'}\, {\delta }_{n,n'}$. Then $(\mathfrak{H}, \langle \, \, , \, \, \rangle )$  is 
a Hilbert space. On $\mathfrak{H}$ we have two self adjoint operators ${\mathbf{Q}}_{A_j}$, $j=1,2$, where 
\begin{equation}
{\mathbf{A}}_1({\mathbf{e}}_{m,n}) = m\hbar \, {\mathbf{e}}_{m,n} \, \, \, \mathrm{and} \, \, \, 
{\mathbf{A}}_2({\mathbf{e}}_{m,n}) = n\hbar \, {\mathbf{e}}_{m,n}.
\label{eq-secn3newtwo}
\end{equation}
The quantum spectrum of the harmonic oscillator is the joint spectrum $({\Z }_{\ge 0})^2$ of 
the operators ${\mathbf{Q}}_{A_j}$, $j=1,2$.  
\vspace{.2in}
\par\noindent \hspace{1in} 
\begin{tabular}{l}
\includegraphics[width=200pt]{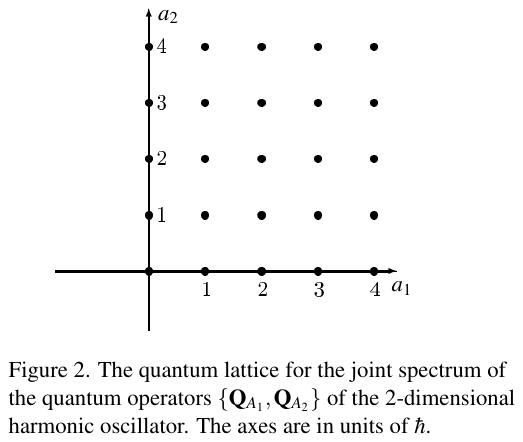} \\ 
\end{tabular} 

To complete the Bohr-Sommerfeld quantization to the Bohr-Sommerfeld-Heisenberg quantization we 
define lowering operators ${\mathbf{a}}_j$ on $\mathfrak{H}$ as follows 
\begin{equation}
{\mathbf{a}}_1{\mathbf{e}}_{m,n} = \left\{ \hspace{-5pt}\begin{array}{cl} {\mathbf{e}}_{m-1,n}, & 
\hspace{-5pt}\mbox{if $m \in {\Z }_{\ge 1}$} \\
0, & \hspace{-5pt}\mbox{if $m=0$} \end{array} \right.  \mathrm{and} \, \, 
{\mathbf{a}}_2{\mathbf{e}}_{m,n} = \left\{ \hspace{-5pt} \begin{array}{cl} {\mathbf{e}}_{m,n-1}, & 
\hspace{-5pt}\mbox{if $n \in {\Z }_{\ge 1}$} \\
0, & \hspace{-5pt}\mbox{if $n=0$.} \end{array} \right. 
\label{eq-secn3newthree}
\end{equation}
According to the theory of section 1, these lowering operators ${\mathbf{a}}_j$, $j=1,2$, correspond to 
the quantum operators ${\mathbf{Q}}_{{\re }^{-\ri \, {\vartheta }_j}}$, $j=1,2$, respectively. Here 
${\vartheta }_j$ is the angle coordinate corresponding to the action $A_j$. However, the functions 
${\re }^{-\ri \, {\vartheta }_j}$ are not smooth on phase space. To remedy this difficiency, we use  
the functions $z_j = {\xi }_j + \ri \, {\eta }_j = r_j{\re }^{-\ri \, {\vartheta }_j}$, $j=1,2$, instead. Then 
${\mathbf{Q}}_{z_j}$, $j=1,2$, are lowering operators on $\mathfrak{H}$, where 
\begin{subequations}
\begin{equation}
{\mathbf{Q}}_{z_1}{\mathbf{e}}_{m,n} = \left\{ \begin{array}{rl} b^1_{m,n} \, {\mathbf{e}}_{m-1,n}, & 
\mbox{if $m \in {\Z }_{\ge 1}$} \\ 0, & \mbox{if $m=0$} \end{array} \right. ,
\label{eq-secn3newfoura}
\end{equation}
that is, $b^1_{0,n} =0$, and
\begin{equation} 
{\mathbf{Q}}_{z_2}{\mathbf{e}}_{m,n} = \left\{ \begin{array}{rl} b^2_{m,n} \, {\mathbf{e}}_{m,n-1}, & 
\mbox{if $n \in {\Z }_{\ge 1}$} \\ 0, & \mbox{if $n=0$,} \end{array} \right. , 
\label{eq-secn3newfourb}
\end{equation}
\end{subequations}
that is, $b^2_{m,0} =0$. The $b^j_{m,n} \in \mathbb{C}$ are determined below. On $\mathfrak{H}$ we have the quantum operator ${\mathbf{Q}}^{\dagger}_{z_j}$, $j=1,2$, which is the adjoint of ${\mathbf{Q}}_{z_j}$. 
By definition ${\mathbf{Q}}^{\dagger}_{z_1} {\mathbf{e}}_{m,n} = {\overline{b}}^1_{m,n} {\mathbf{e}}_{m+1,n}$ 
and ${\mathbf{Q}}^{\dagger}_{z_2} {\mathbf{e}}_{m,n} = {\overline{b}}^2_{m,n} {\mathbf{e}}_{m, n+1}$. For 
$j=1,2$ we identify ${\mathbf{Q}}^{\dagger}_{z_j}$ with ${\mathbf{Q}}_{{\overline{z}}_j}$. \medskip 

For $1 \le j,k \le 2$ we have 
\begin{displaymath}
\{ z_j, {\overline{z}}_k \} = \{ {\xi }_j + \ri \, {\eta }_j, {\xi }_k - \ri \, {\eta }_k \} = 
-\ri \, \{ {\xi }_j, {\eta }_k \} + \ri \, \{ {\eta }_k, {\xi }_j \} = -2\ri \, {\delta }_{j,k} 
\end{displaymath}
and $\{ z_j, z_k \} = 0 = \{ {\overline{z}}_j, {\overline{z}}_k \} $. So by the Dirac quantization prescription we get 
\begin{subequations}
\begin{align}
[{\mathbf{Q}}_{z_j}, {\mathbf{Q}}_{A_k} ] & = [ {\mathbf{Q}}_{z_j}, {\mathbf{Q}}_{\ttfrac{1}{4}  (z_k {\overline{z}}_k  + {\overline{z}}_k z_k)}] = \ri \hbar \, {\mathbf{Q}}_{ \ttfrac{1}{4} \{ z_j, z_k {\overline{z}}_k +{\overline{z}}_k z_k \} } \notag \\
& = \ri \hbar \, {\mathbf{Q}}_{\onehalf z_k \{ z_j, {\overline{z}}_k \} } 
= \ri \hbar \, {\mathbf{Q}}_{-\ri \, z_k \, {\delta }_{j,k} } = 
\hbar \, {\delta }_{j,k} \, {\mathbf{Q}}_{z_k}; 
\label{eq-secn3newfivea} 
\end{align}
\begin{align}
[{\mathbf{Q}}_{{\overline{z}}_j}, {\mathbf{Q}}_{A_k} ] & = \ri \hbar \, {\mathbf{Q}}_{\onehalf {\overline{z}}_k \{ {\overline{z}}_j, z_k \} } 
= \ri \hbar \, {\mathbf{Q}}_{\ri \, {\overline{z}}_k \, {\delta }_{j,k} } = 
-\hbar \, {\delta }_{j,k} \, {\mathbf{Q}}_{{\overline{z}}_k};  
\label{eq-secn3newfiveb}  
\end{align}
and 
\begin{align}
[{\mathbf{Q}}_{z_j}, {\mathbf{Q}}_{{\overline{z}}_k}] & = \ri \hbar \, {\mathbf{Q}}_{\{ z_j, {\overline{z}}_k \} } 
= \ri \hbar \, {\delta }_{j,k} \, {\mathbf{Q}}_{-2\ri \, {\delta }_{j,k}} = 2\hbar \, {\delta  }_{j,k}.  
\label{eq-secn3newfivec}  
\end{align}
\end{subequations}
Thus $\{ {\mathbf{Q}}_{z_j}, {\mathbf{Q}}_{{\overline{z}}_j}, {\mathbf{Q}}_{A_j} \} $, $j=1,2$, span 
a finite dimensional Lie algebra with bracket relations 
\begin{displaymath}
[{\mathbf{Q}}_{z_j}, {\mathbf{Q}}_{A_j} ] = 
\ri \hbar \, {\mathbf{Q}}_{z_j}, \, \,  [{\mathbf{Q}}_{{\overline{z}}_j}, {\mathbf{Q}}_{A_j} ] = 
-\ri \hbar \, {\mathbf{Q}}_{{\overline{z}}_j}, \, \, \mathrm{and} \, \,  [ {\mathbf{Q}}_{z_j}, {\mathbf{Q}}_{{\overline{z}}_j} ] = 2\hbar , 
\end{displaymath}
which commute with each other. \medskip 

The quantum operators ${\mathbf{Q}}_{z_j}$ 
and ${\mathbf{Q}}_{{\overline{z}}_j}$ are not self adjoint. However, we can define  
the quantum operators ${\mathbf{Q}}_{{\xi}_j} =
\onehalf ({\mathbf{Q}}_{z_j} + {\mathbf{Q}}_{{\overline{z}}_j})$ and 
${\mathbf{Q}}_{{\eta }_j} = \ttfrac{1}{2 \ri} ({\mathbf{Q}}_{z_j} - {\mathbf{Q}}_{{\overline{z}}_j})$, 
which are self adjoint. A straightforward calculation shows that the skew hermitian operators 
$\{  \ttfrac{1}{\ri \hbar}\, {\mathbf{Q}}_{{\xi }_j}, \, 
\ttfrac{1}{\ri \hbar}\, {\mathbf{Q}}_{{\eta }_j}, \, 
\ttfrac{1}{\ri \hbar}\, {\mathbf{Q}}_{A_j} \} $ satisfies the bracket relations 
\begin{subequations}
\begin{align}
[ \ttfrac{1}{\ri \hbar} \, {\mathbf{Q}}_{{\xi }_j}, 
\ttfrac{1}{\ri \hbar}\, {\mathbf{Q}}_{A_j} ] & =
\ttfrac{1}{\ri \hbar} \, {\mathbf{Q}}_{{\eta }_j},
\label{eq-secn3newsixa} \\
\lbrack \ttfrac{1}{\ri \hbar}\,  {\mathbf{Q}}_{ {\eta }_j}, 
\ttfrac{1}{\ri \hbar}\, {\mathbf{Q}}_{A_j} \rbrack & = 
- \ttfrac{1}{\ri \hbar}\, {\mathbf{Q}}_{ {\xi }_j} ,
\label{eq-secn3newsixb} \\  
[ \ttfrac{1}{\ri \hbar}\, {\mathbf{Q}}_{{\xi }_j}, 
\ttfrac{1}{\ri \hbar}\, {\mathbf{Q}}_{{\eta }_j} ] & = 
\ttfrac{1}{\ri \hbar} 
\label{eq-secn3newsixc}
\end{align}
\end{subequations}
and thus span a finite dimensional Lie algebra ${\mathcal{L}}_j$, $j=1,2$, which is a subalgebra of 
$\mathfrak{u}(\mathfrak{H}, \mathbb{C})$, and $[{\mathcal{L}}_1, {\mathcal{L}}_2] =0$. 
Note that $2{\mathbf{Q}}_{z_j} = {\mathbf{Q}}_{{\xi }_j} + \ri \, {\mathbf{Q}}_{{\eta }_j}$ and 
$2\ri \, {\mathbf{Q}}_{{\overline{z}}_j} = {\mathbf{Q}}_{{\xi }_j} - \ri \, {\mathbf{Q}}_{{\eta }_j}$.  
\medskip 

We now determine $b^j_{m,n}$ in (\ref{eq-secn3newfoura}) and (\ref{eq-secn3newfourb}). Suppose 
for every $(m,n) \in ({\Z }_{\ge 0})^2$ that 
$b^j_{m,n} \in \R $. When $m \in {\Z }_{\ge 1}$ we have 
\begin{align}
{\mathbf{Q}}_{z_1{\overline{z}}_1} {\mathbf{e}}_{m,n} & = 
{\mathbf{Q}}_{z_1}{\mathbf{Q}}_{{\overline{z}}_1} {\mathbf{e}}_{m,n} = 
b^1_{m+1,n}\, {\mathbf{Q}}_{z_1} {\mathbf{e}}_{m+1,n} = (b^1_{m+1,n})^2 \, {\mathbf{e}}_{m,n} \notag \\
{\mathbf{Q}}_{{\overline{z}}_1z_1} {\mathbf{e}}_{m,n} & = (b^1_{m,n})^2 \, {\mathbf{e}}_{m,n} \notag \\ 
{\mathbf{Q}}_{z_1{\overline{z}}_1} {\mathbf{e}}_{m,n} & = 
{\mathbf{Q}}_{z_2}{\mathbf{Q}}_{{\overline{z}}_2} {\mathbf{e}}_{m,n} = 
(b^2_{m,n+1})^2 \, {\mathbf{e}}_{m,n} \notag \\
{\mathbf{Q}}_{{\overline{z}}_2z_2} {\mathbf{e}}_{m,n} & = (b^2_{m,n})^2 \, {\mathbf{e}}_{m,n}. \notag 
\end{align}
Therefore we get 
\begin{align}
4m\hbar \, {\mathbf{e}}_{m,n} & = {\mathbf{Q}}_{4A_1}{\mathbf{e}}_{m,n} = {\mathbf{Q}}_{z_1{\overline{z}}_1} {\mathbf{e}}_{m,n} + {\mathbf{Q}}_{{\overline{z}}_1z_1} {\mathbf{e}}_{m,n} \notag \\  
& = \big( (b^1_{m+1,n})^2 +(b^1_{m,n})^2 \big) {\mathbf{e}}_{m,n} , \notag 
\end{align}
which implies \begin{subequations}
\begin{equation}
4m \hbar =  (b^1_{m+1,n})^2 +(b^1_{m,n})^2  , 
\label{eq-secn3newsevena}
\end{equation}
and
\begin{align}
2\hbar \, {\mathbf{e}}_{m,n} & = [ {\mathbf{Q}}_{z_1}, {\mathbf{Q}}_{{\overline{z}}_1} ] {\mathbf{e}}_{m,n} 
= {\mathbf{Q}}_{z_1{\overline{z}}_1} {\mathbf{e}}_{m,n} - 
{\mathbf{Q}}_{{\overline{z}}_1z_1} {\mathbf{e}}_{m,n} \notag \\ 
& = \big( (b^1_{m+1,n})^2 - (b^1_{m,n})^2 \big) {\mathbf{e}}_{m,n}, \notag 
\end{align}
which implies 
\begin{equation}
2\hbar = (b^1_{m+1,n})^2 - (b^1_{m,n})^2. 
\label{eq-secn3newsevenb}
\end{equation}
\end{subequations}
Subtracting (\ref{eq-secn3newsevenb}) from (\ref{eq-secn3newsevena}) gives $b^1_{m,n} = \sqrt{(2m-1)\hbar }$, if 
$m \in {\Z }_{\ge 1}$. Recall that by definition $b^1_{0,n} =0$. Arguing similarly we get $b^2_{m,n}= \sqrt{(2n-1)\hbar }$ if $n \in {\Z }_{\ge 1}$ and $b^2_{m,0} =0$ by definition. 
This completes the definition of the shifting operators on $\mathfrak{H}$. \medskip 

We now look at the quantum representation of the Lie algebra ${\mathcal{L}}_2$ on the Hilbert space 
$(\mathfrak{H}, \langle \, \, , \, \, \rangle )$ defined by the skew hermitian 
linear maps $\{  \tttfrac{1}{\ri \hbar} \, {\mathbf{Q}}_{{\mathrm{Re}}\, z_2}, $
$ \, \tttfrac{1}{\ri \hbar} \, {\mathbf{Q}}_{{\mathrm{Im}}\, z_2},$ 
$ \, \tttfrac{1}{\ri \hbar} \, {\mathbf{Q}}_{A_2} \} $ of $\mathfrak{H}$ into 
itself and their brackets, which satisfy (\ref{eq-secn3newsixa})--(\ref{eq-secn3newsixc}). For each $m 
\in {\Z }_{\ge 0}$ let  ${\mathfrak{H}}_m = \mathrm{span}\, \{ {\mathbf{e}}_{m,n} \, | \, 
n \in {\Z }_{\ge 0} \} $. Then ${\mathfrak{H}}_m = \mathrm{span}\, \{ 
Q^n_{{\overline{z}}_2}{\mathbf{e}}_{m,0} \, | \, n \in {\mathbf{Z}}_{\ge 0} \} $, since 
$Q^n_{{\overline{z}}_2}{\mathbf{e}}_{m,0} \in \mathrm{span}\, \{ {\mathbf{e}}_{m,n} \}$. Because 
${\mathcal{L}}_2({\mathfrak{H}}_m)$ $\subseteq 
{\mathfrak{H}}_m $, we have a representation of ${\mathcal{L}}_2$ on ${\mathfrak{H}}_m$, which  
is \emph{irreducible} since the vectors ${\{ Q^n_{{\overline{z}}_2}{\mathbf{e}}_{m,0} \} }^{\infty}_{n=0}$ 
span ${\mathfrak{H}}_m $. This representation is a faithful and infinite dimensional with a unique lowest weight vector ${\mathbf{e}}_{m,0}$, but no top weight vector. Since $\mathfrak{H} = 
\bigoplus_{m \in {\mathbf{Z}}_{\ge 0}} {\mathfrak{H}}_m $, the representation of ${\mathcal{L}}_2$ on 
$\mathfrak{H}$ is a countably infinite number of copies of the irreducible representation of 
${\mathcal{L}}_2$ on ${\mathfrak{H}}_0$. A similar argument can be made for the quantum representation of 
${\mathcal{L}}_1$ on $\mathfrak{H}$. 

 \section{The reduction of the oscillator symmetry}

In this section we reduce the $S^1$ symmetry of the harmonic oscillator generated by the flow of 
the harmonic oscillator vector field on phase space. 
\medskip

For $j=1,2$ let $z_{j}={\xi }_{j}+ \ri \,{\eta }_{j} $. Then 
${\R }^{4}$ becomes ${\mathbb{C}}^{2}$ with coordinates $(z_{1},z_{2})$. The flow 
${\varphi }_{t}^{E}$ of the vector field $X_{E}$ becomes the complex $S^{1}$-action 
\begin{equation*}
\varphi :S^{1}\times {\mathbb{C}}^{2}\rightarrow {\mathbb{C}}^{2}:
\big( t,(z_{1},z_{2})\big)\mapsto ({\re }^{-\ri t}z_{1},{\re }^{-\ri t}z_{2}).
\end{equation*}

The algebra of $S^{1}$-invariant polynomials on ${\R }^{4}$ is
generated by 
\begin{align}
{\pi }_1 & = \mathrm{Re}({\overline{z}}_1z_2) = {\xi }_1{\xi }_2+{\eta }_1{\eta }_2 \notag \\    
{\pi }_2 & = \mathrm{Im}({\overline{z}}_1z_2) =  {\xi }_1{\eta }_2-{\xi }_2{\eta }_1 
\label{eq-sec1invariants} \\  
{\pi }_3 & = \onehalf  (z_1{\overline{z}}_1-z_2{\overline{z}}_2)  = 
\onehalf  ({\xi }^2_1+{\eta }^2_1-{\xi }^2_2 - {\eta }^2_2) \notag   \\  
{\pi }_4 & = \onehalf  (z_1{\overline{z}}_1 +z_2{\overline{z}}_1) = 
\onehalf  ({\xi }^2_1+{\eta }^2_1+{\xi }^2_2 + {\eta }^2_2) \notag 
\end{align}%
subject to the relation 
\begin{equation*}
{\pi }_{1}^{2}+{\pi }_{2}^{2}=|{\overline{z}}_{1}z_{2}|^{2}
=(z_{1}{\overline{z}}_{1})(z_{2}{\overline{z}}_{2})
=({\pi }_{4}+{\pi }_{3})({\pi }_{4}-{\pi }_{3})
={\pi }_{4}^{2}-{\pi }_{3}^{2},\mspace{10mu} {\pi }_{4}\geq 0,
\end{equation*}
which defines the orbit space ${\R }/S^1$. Note that the energy $E$ (\ref{eq-sec2newzeroa}) is ${\pi }_4$ and 
the angular momentum $L$ (\ref{eq-sec2newzerob}) is ${\pi }_3$.
\medskip

Using the Poisson bracket on $C^{\infty}({\R  }^4)$ we find the
structure matrix $(\{ {\pi}_i, {\pi}_j \})$ of the Poisson bracket
on ${\R  }^4$ with coordinates $({\pi }_1, \ldots , {\pi }_4)$,
namely, \bigskip

\noindent \hspace{1.25in}\begin{tabular}{c|cccc}
$\{{\pi }_i,{\pi }_j\}$ & ${\pi }_1$ & ${\pi }_2$ & ${\pi }_3$ & ${\pi }_4$
\\ \hline
${\pi }_1$ & $0$ & $2{\pi }_3$ & $-2{\pi }_2$ & $0$ \\ 
${\pi }_2$ & $-2{\pi }_3$ & $0$ & $2{\pi }_1$ & $0$ \\ 
${\pi }_3$ & $2{\pi }_2$ & $-2{\pi }_1$ & $0$ & $0$ \\ 
${\pi }_4$ & $0$ & $0$ & $0$ & $0$
\end{tabular}
\vspace{.2in}

\noindent \hspace{1.2in}\parbox[t]{4in}{Table 1. The Poisson bracket on 
$C^{\infty}({\R }^4)$.}\medskip

Let $S_{\sqrt{2e}}^3 =\{ (\xi , \eta ) \in {\R }^4 \, | \, {\xi}^2_1 + {\eta }^2_1+ {\xi }^2_2 +{\eta }^2_2 =2e \}$ 
be the $3$-sphere of radius $\sqrt{2e}$ and let $S^2_e = 
\{ \pi \in {\R  }^3 \, | \, {\pi }^2_1+{\pi }^2_2+ {\pi }^2_3 =e^2\}$ be the $2$-sphere of
radius $e$. The map 
\begin{equation}
\rho :S_{\sqrt{2e}}^3 \subseteq {\R  }^4\rightarrow S^2_e \subseteq {\R  }^3 : 
\zeta =(\xi ,\eta )\mapsto \pi = ({\pi }_1(\zeta ), {\pi }_2(\zeta ),{\pi }_3(\zeta ))  
\label{hopf}
\end{equation}
is called the Hopf fibration, see \cite[chapter I]{cushman-bates}. 
Each fiber of the Hopf fibration is an integral curve of the harmonic
oscillator vector field $X_E$ of energy $e$. Thus the space $E^{-1}(e)/S^1$
of orbits of the harmonic oscillator of energy $e$ is $S_e^2$. In other
words, the Hopf fibration $\rho $ (\ref{hopf}) is the reduction map of the $S^1$ symmetry of the harmonic oscillator generated by the vector field $X_E$. \medskip

Consider a smooth function $K:{\R  }^4 \to \R $, which is
invariant under the flow of $X_E$. Then $K$ is an integral of $X_E$, that
is, ${\pounds }_{X_E}K=0$. Also, there is a smooth function $\widetilde{K}:
{\R }^4 \to \R $ such that $K(\xi , \eta )= 
\widetilde{K}({\pi}_1(\xi ,\eta ),{\pi }_2(\xi ,\eta ),$ ${\pi }_3(\xi ,\eta ),$ ${\pi }_4(\xi ,\eta ))$. \medskip

Since $\{{\pi }_j,{\pi}_4\}=0$ for $j=1,\dotsc ,4$, we obtain 
\begin{equation*}
\dot{\pi }_j = \{ {\pi}_j, \widetilde{K} \} \, = \, \sum_{k=1}^{3}\{{\pi }_j,{\pi }_k\} 
\frac{\partial\widetilde{K}}{\partial {\pi }_k} = 2
\sum_{k=1}^{3}\sum_{l=1}^{3} \varepsilon_{jkl}\frac{\partial\widetilde{K}}{\partial {\pi }_k}\, {\pi}_l 
= 2(\nabla \widetilde{K}\times \pi )_j
\end{equation*}
for $j=1,2,3$ and $\dot{\pi}_4=0$. \medskip

Restrict the function $\widetilde{K}$ to $E^{-1}(e)$ and define $\widetilde{K}_e({\pi }_1,{\pi }_2, {\pi }_3)$ 
$= \widetilde{K}({\pi}_1,{\pi }_2,{\pi }_3, $ $e)$. Set $\pi =({\pi }_1,{\pi }_2,{\pi} _3) \in {\R  }^3$ and let 
$\langle \, \, , \, \, \rangle $ be the Euclidean inner product on ${\R }^3$. Then $\dot{\pi} = 2(\nabla{\widetilde{K}}_e\times \pi )$ is satisfied by
integral curves of a vector field $X$ on ${\R }^3$ defined by 
\begin{equation}
X(\pi ) = 2 (\nabla{\widetilde{K}}_e \times \pi ).  
\label{eq-hamvf}
\end{equation}
The $2$-sphere $S_e^2$ is invariant under the flow of the vector field $X$,
because 
\begin{equation*}
{\pounds }_{X}\langle \pi , \pi \rangle = 2\langle \pi, \dot{\pi}\rangle =
4\langle \pi , \nabla{\widetilde{K}}_e(\pi )\times \pi \rangle =0.
\end{equation*}

The structure matrix of the Poisson bracket $\{\,\,,\,\,\}$ on $C^{\infty }({\R }^{3})$ is 
\begin{equation}
W(\pi )=(\{{\pi }_{j},{\pi }_{k}\})=
-2\mbox{\footnotesize $\begin{pmatrix}
 0 & -{\pi }_3 & {\pi }_2 \\ {\pi }_3 & 0 & -{\pi }_1 \\ -{\pi }_2 & {\pi }_1 & 0
 \end{pmatrix}$}
\end{equation}%
Since $\ker W(\pi )=\mathrm{span} \{\pi \}$ and $T_{\pi}S_{e}^{2}={\mathrm{span} \{\pi \}}^{\perp }$, the matrix $W(\pi )|_{T_{\pi }S_{e}^{2}}$ is invertible. On $S_{e}^{2}$ define a symplectic
form ${\widetilde{\omega }}_{e}({\pi })(u,v)=\langle (W({\pi })^{T})^{-1}u,v\rangle $,
where $u,v\in T_{\pi }S_{e}^{2}$. Let $y\in T_{\pi }S_{e}^{2}$. Since $W(\pi )^{T}y=2{\pi }\times {}y=z$ we get 
\begin{align}
\pi \times {}z& =\pi \times (2\pi \times {}y)=2\pi \times (\pi \times {}y) 
\notag \\
& =2(\pi \,\langle \pi ,y\rangle -y\,\langle \pi ,\pi \rangle )=-2y\,\langle
\pi ,\pi \rangle =-2e^{2}\,y,  \notag
\end{align}%
which implies $y=(W(\pi )^{T})^{-1}z=- \tfrac{1}{2e^2} \, \pi \times {}z$. 
Therefore 
\begin{equation}
{\widetilde{\omega}} _{e}(\pi )(u,v)=
-\tfrac{1}{2e^2} \, \langle \pi \times {}u,v\rangle =
-\tfrac{1}{2e^2} \, \langle \pi ,u\times {}v\rangle .
\end{equation}

The vector field $X$ (\ref{eq-hamvf}) on $(S^2_e, {\widetilde{\omega}}_e)$ is Hamiltonian
with Hamiltonian function $\widetilde{K}_e$, because 
\begin{align}
{\widetilde{\omega}}_e(\pi )(X(\pi ),u) & = - \tfrac{1}{2e^2} \, 
\langle \pi , 2(\nabla \widetilde{K}_e\times \pi ) \times u \rangle
= - \tfrac{1}{e^2} \, \langle \pi \times
(\nabla \widetilde{K}_e\times \pi ), u \rangle  \notag \\
& = - \tfrac{1}{e^2} \, \langle \nabla 
\widetilde{K}_e \langle \pi , \pi \rangle - \pi \langle \nabla {\widetilde{K}}_e, \pi \rangle , u \rangle 
= -\langle \nabla\widetilde{K}_e, u\rangle \notag \\
& = -\mathrm{d} \widetilde{K}_e(\pi ) u,  \notag
\end{align}
where $u,v\in T_w S_e^2$. \medskip

Since $L$ is a function which is invariant under the flow ${\varphi }^E_t$
of $X_E $ on $E^{-1}(e)$, it induces the function ${\widetilde{L}}_e = {\widetilde{\pi }}_3 = {\pi }_3|S^2_e$ on the 
reduced phase space $S^2_e$. We look at the completely
integrable reduced system $({\widetilde{L}}_e, S^2_e, {\widetilde{\omega }}_e)$.
\medskip

The flow of the reduced vector field $X_{{\widetilde{L}}_e}(\pi ) = 2
(e_3 \times \pi )$ on $S^2_e \subseteq {\R  }^3$ is 
\begin{equation*}
{\varphi }^{{\widetilde{\pi }}_3}_t(\pi ) =  
\mbox{{\footnotesize $\begin{pmatrix} \cos 2t &- \sin 2t & 0 \\ \sin 2t & \cos 2t & 0 \\
0 & 0 & 1 \end{pmatrix}$}}\pi ,
\end{equation*}
where $\pi \in S^2_e$. Note that ${\varphi }^{{\widetilde{L}}_e}_t$ is periodic of period $\pi $. 
Therefore an integral curve of $X_{{\widetilde{\pi }}_3} $ on $S^2_e$
is either a circle or a point. For $|\ell | <e$ we have ${\pi }^{-1}_3(\ell ) $ is the circle 
$\{ ({\pi }_1, {\pi }_2, \ell ) \in S^2_e \, | \, {\pi }^2_1+{\pi }^2_2 = e^2- {\ell }^2 \} $; while 
when $\ell = \pm e$, we see that ${\pi}^{-1}(\ell )$ is the point $(0,0,\pm e)$. \medskip

Introduce spherical coordinates 
\begin{equation*}
{\pi }_{1}=e\,\sin \theta \cos 2\psi ,\quad {\pi }_{2}=e\sin \theta \sin
2\psi ,\mspace{5mu}\text{and}\mspace{5mu} {\pi }_{3}=e\cos \theta
\end{equation*}
with $0\leq \theta \leq \pi $ and $0\leq \psi \leq \pi $. Therefore the
reduced symplectic form on $S_{e}^{2}$ is 
\begin{equation*}
{\widetilde{\omega }}_{e}=-\ttfrac{1}{2e}\, (e^{2}\sin \theta )\rd  
\theta \wedge \rd  \psi ,
\end{equation*}
that is, ${\widetilde{\omega }}_{e}=\ttfrac{1}{2e} \, {\mathrm{vol}}_{S_{e}^{2}}$. Here ${\mathrm{vol}}_{S_{e}^{2}}$ is the standard volume $2$-form $- \ttfrac{1}{2e} \, \langle \pi ,u\times v\rangle $ on $S_{e}^{2}$, where 
$\pi \in S_{e}^{2}$
and $u,v\in T_{\pi }S_{e}^{2}$. Note $\int_{S_{e}^{2}}{\mathrm{vol}}_{S_{e}^{2}}
=4\pi e^{2}$. \medskip

Following our convention, for $i=1,2,3,$ we denote by $\widetilde{\pi}_{i}$
the push-forward to $S_{e}^{2}$ of the invariant function $\pi _{i}$ on $%
\R ^{3}$. We now explain why the functions ${\widetilde{\pi }}_i$
satisfy commutation relations for $\su (2)$.
Consider the Poisson \linebreak 
algebra $\mathcal{A} = (C^{\infty}({\R  }^3), {%
\{ \, \, , \, \, \}}_{{\R  }^3}, \cdot )$, where ${\R  }^3$
has coordinates $({\pi }_1, {\pi }_2, {\pi }_3)$, the Poisson bracket 
${\{ \, \, , \, \, \}}_{{\R  }^3}$ has structure matrix 
$(\{{\pi }_{j},{\pi }_{k}\}) =${\tiny $%
\begin{pmatrix}
0 & 2{\pi }_3 & -2{\pi }_2 \\ 
-2{\pi }_3 & 0 & 2{\pi }_1 \\ 
2{\pi }_2 & -2{\pi }_1 & 0%
\end{pmatrix}%
$}, and $\cdot $ is pointwise multiplication of functions. Since $C = {\pi }^2_1 +{\pi }^2_2+{\pi }^2_3-e^2$ 
is a Casimir of $\mathcal{A}$, that is, ${ \{ C , {\pi }_i \} } _{{\R  }^3} =0$ for $i=1,2,3$,
it generates a Poisson ideal $\mathcal{I}$ of $\mathcal{A}$. Observe that $\mathcal{I} = 
\{ f \in C^{\infty}({\R  }^3) |\,  f |_{S^2_e} =0 \} $, where $S^2_e =
C^{-1}(0)$. Therefore $\mathcal{B} = (C^{\infty}(S^2_e) = C^{\infty}({%
\R  }^3)/\mathcal{I}, {\{ \, \, , \, \, \} }_{S^2_e}, \cdot )$ is a
Poisson  
algebra, where ${ \{ \overline{f} , \overline{g} \} }_{S^2_e} = \mbox{${\{ f , g \}} _{{\scriptstyle {\R  }^3}}$}%
|_{S^2_e}$ and $f \mapsto \overline{f}$ is the natural projection map of $C^{\infty}({\R  }^3)$ 
onto $C^{\infty}(S^2_e)$. Note that $\overline{f} = f|_{S^2_e}$, because 
every smooth function on ${\R  }^3$, which
vanishes on $S^2_e$, is a multiple of $C$ by a smooth function on ${\R  }^3$. Consequently, the structure matrix of the Poisson bracket ${\{ \, \, , \, \, \}}_{S^2_e}$ is {\tiny $%
\begin{pmatrix}
0 & 2{\widetilde{\pi }}_3 & -2{\widetilde{\pi }}_2 \\ 
-2{\widetilde{\pi }}_3 & 0 & 2{\widetilde{\pi }}_1 \\ 
2{\widetilde{\pi }}_2 & -2{\widetilde{\pi }}_1 & 0%
\end{pmatrix}%
$}, where ${\widetilde{\pi }}_i = {\pi }_i|_{S^2_e}$. In other words, ${\{ {%
\widetilde{\pi }}_i, {\widetilde{\pi }}_j \} }_{S^2_e} = 2\sum^3_{k=1}{%
\varepsilon }_{ijk}{\widetilde{\pi }}_k$, which are the bracket relations
for $\su (2)$, as desired. The reduced phase space $S^2_e$ is symplectomorphic to a
coadjoint orbit of $\mathrm{SU}(2)$.  \medskip

The Hamiltonian vector field $X_{{\widetilde{L}}_e}$ is $\frac{\partial }{\partial \psi }$, because 
\begin{equation*}
X_{{\widetilde{L}}_e}\lefthook \, {\widetilde{\omega }}_{e} 
=e\sin \theta \, \rd  \theta =-\rd  (e\cos \theta )=
-\rd  {\widetilde{\pi }}_{3}.
\end{equation*}%
Note that the flow of $X_{{\widetilde{L}}_e}$ on $S_{e}^{2}$ is periodic of
period $\pi $. On $S_{e}^{2}\setminus \{(0,0,\pm e)\}$ we have 
\begin{equation*}
\rd  {\widetilde{L}}_e \wedge \rd  \psi = \rd  {\widetilde{\pi }}_{3}\wedge \rd  \psi =
\rd  (e\cos \theta )\wedge \rd  \psi = 
-e\sin \theta \, \rd  \theta \wedge \rd  \psi  ={\widetilde{\omega }}_{e}.
\end{equation*}%
Thus $({\widetilde{L}}_e, \psi )$ are real analytic action-angle coordinates
on $S_{e}^{2}\setminus \{(0,0,\pm e)\} $ with the symplectic form ${\widetilde{\omega }}_e|_{S_{e}^{2}\setminus 
\{ (0,0,\pm e) \} }$.

\section{Quantization of the reduced system}

For fixed $e\geq 0$ we quantize the classical system $({\widetilde{L}}_e,S_{e}^{2},{\widetilde{\omega }}_{e})$ 
obtained by reducing the harmonic oscillator symmetry of the Hamiltonian system $(L,T^{\ast }{\R }^{2},$ $\omega )$ on $E^{-1}(e)$, see \S 3. We find that 
its space of quantum states is the sum of two finite dimensional irreducible $\su (2)$ representations: 
one of dimension $q$ and the other of dimension $q+1$. The removal of all infinite dimensional representations and most of their multiplicities is mainly due to restricting the harmonic oscillator to a fixed energy level. \medskip

Recall that $(S_{e}^{2},{\widetilde{\omega }}_{e})$ is symplectomorphic to a coadjoint
orbit of $\mathrm{SU}(2)$. Moreover, functions $\widetilde{\pi}_{1},\widetilde{\pi}_{2},\widetilde{\pi}_{3}$ 
correspond to evaluation of elements of the orbit on a basis of $\su (2)$. Hence, quantization of the system 
$({\widetilde{L}}_e = {\widetilde{\pi }}_{3},S_{e}^{2},{\widetilde{\omega }}_{e})$ should lead to a
representation of $\mathrm{SU} (2).$\medskip

We assume that the symplectic manifold $(S_{e}^{2}, {\widetilde{\omega }}_{e})$ is
prequantizable. In other words, for $q\in \Z $ we have 
\begin{equation*}
qh=\int_{S_{e}^{2}}{\widetilde{\omega }}_{e}= \ttfrac{1}{2e} \, \int_{S_{e}^{2}}{\mathrm{vol}}_{S_{e}^{2}}= 
\ttfrac{1}{2e} (4\pi e^{2})=2\pi e,
\end{equation*}%
that is, $0\leq e=q\hbar $, where $\hbar =\ttfrac{h}{2\pi }$. Thus $q$ is a fixed nonnegative integer.
\medskip

The Bohr-Sommerfeld conditions for the integrable reduced system $({\widetilde{L}}_e$, 
$S_{e}^{2},$ ${\widetilde{\omega }}_{e})$ in action-angle
coordinates $({\widetilde{L}}_e ,\psi )$, see \S 3, are for $p\in \Z $ 
\begin{equation*}
ph=\oint_{a^{-1}(\ell )} {\widetilde{L}}_e  \,\rd  \psi =\ell \int_{0}^{\pi }2\, \rd  \psi =2\pi \ell .
\end{equation*}%
The second equality above follows because the curve 
\begin{displaymath}
[0,\pi ]\rightarrow
S_{e}^{2}:\psi \mapsto \big(\sqrt{e^{2}-{\ell }^{2}}\cos 2\psi ,\sqrt{e^{2}-{\ell }^{2}}\sin 2\psi ,\ell \big) 
\end{displaymath} 
parametrizes $a^{-1}(\ell )$. Thus $\ell =p\hbar $. Since $\ell ={\pi }_{3}\cos {\theta }_{p}$ and 
$|\cos {\theta }_{p}|\leq 1$, it follows that $|\ell |\leq e$, which implies $|p|\leq q$ and 
$\cos {\theta }_{p}=\ttfrac{p}{q}$. So the Bohr-Sommerfeld set 
${\widetilde{S}}_{q}$ for the reduced integrable system 
$({\pi }_{3}|S_{q\hbar }^{2},S_{q\hbar }^{2}, {\widetilde{\omega }}_{q\hbar })$ with $q\in {\Z }_{\geq 0}$ is 
\begin{equation}
{\widetilde{S}}_{q}=\{\big(q\hbar \sin {\theta }_{p}\cos 2\psi ,q\hbar \sin {\theta }_{p}\sin 2\psi ,p\hbar \big)\in 
S_{q\hbar }^{2} \, | \, p \in \Z  \} 
\label{eq-bscircle}
\end{equation}%
with $|p|\leq q$.
The set ${\widetilde{S}}_{q}$ is a disjoint union of $2q+1$ circles of
radius $(\sqrt{{q}^{2}-{p}^{2}})\hbar $ when $|p|<q$ and two points $%
(0,0,\pm q\hbar )$ when $p=\pm q$. \medskip

Let ${\widetilde{\mathbf{e}}}_{p,q}$ with $|p|\leq q$ be a vector in ${\widetilde{\mathfrak{H}}}_{q}$ which corresponds to the $1$-torus 
\begin{displaymath}
\{ \big(q\hbar \sin {\theta }_{p}\cos 2\psi ,q\hbar \sin {\theta }_{p}\sin 2\psi ,p\hbar \big) 
\in S^2_{q\hbar } \, | \, \psi \in [0, \pi ] \}
\end{displaymath} 
in the Bohr-Sommerfeld set ${\widetilde{S}}_{q}$. Define an
inner product $\langle \,\, ,\,\,\rangle $ on ${\widetilde{\mathfrak{H}}}_{q}$
so that $\langle {\widetilde{\mathbf{e}}}_{p,q}, 
{\widetilde{\mathbf{e}}}_{p^{\prime},q^{\prime }}\rangle ={\delta }_{(p,q),(p^{\prime },q^{\prime })}$. 
From now on we use the notation ${\widetilde{\pi}}_{j}$, $j=1,2,3$ for 
${\pi }_{j}|S_{q\hbar }^{2}$. For each integer $p$ with $|p|\leq q$ we have 
\begin{equation}
{\mathbf{Q}}_{{\widetilde{\pi}}_{3}}\widetilde{\mathbf{e}}_{p,q}= q\hbar \cos {\theta 
}_{p}\,\widetilde{\mathbf{e}}_{p,q}= p\hbar \,\widetilde{\mathbf{e}}_{p,q}.
\label{eq-s8twostar}
\end{equation}%
So ${\widetilde{\mathfrak{H}}}_{q} =\mathrm{span} \{{\widetilde{\mathbf{e}}}_{p,q}\,
| \,|p|\leq q\}$ consists of eigenvectors of ${\mathbf{Q}}_{{\widetilde{\pi}}_{3}}$ corresponding 
to the eigenvalue $\hbar p$. \medskip

We want to define shifting operators ${\mathbf{\widetilde{a}}}_{q}$ and 
${\mathbf{\widetilde{a}}}_{q}^{\dagger }$ on 
$\widetilde{\mathfrak{H}}$ so that 
\begin{equation*}
{\mathbf{\widetilde{a}}}_{q}\,\widetilde{\mathbf{e}}_{p,q}=\widetilde{\mathbf{e}}_{p-1,q}%
\mspace{10mu}\text{and}\mspace{10mu} 
{\mathbf{\widetilde{a}}}_{q}^{\dagger}\,{\widetilde{\mathbf{e}}}_{p,q}=
{\widetilde{\mathbf{e}}}_{p+1,q}
\end{equation*}%
By general theory, we expect to identify the operators ${\mathbf{\widetilde{a}}}_{q}$ and 
${\mathbf{\widetilde{a}}}_{q}^{\dagger }$ with the quantum
operators ${\mathbf{Q}}_{{\re }^{-\ri \psi }}$ and ${\mathbf{Q}}_{{\re }^{\ri \psi }}$, respectively. However, the functions 
${\re }^{\pm \ri \psi }$ are not single-valued on the complement of the poles in $S_{q\hbar}^{2}$. In order 
to get single-valued functions, we have take squares of ${\re }^{\pm \ri\, \psi }$, namely the functions 
${\re }^{\pm 2i\, \psi }$. The functions ${\re }^{\pm 2\ri \, \psi }$ do not extend to smooth functions on 
$S_{q\hbar }^{2}$. However, the functions 
\begin{equation*}
{\widetilde{\pi }}_{\pm }= \big( \sqrt{(q\hbar )^{2}-{\widetilde{\pi }}_{3}^{2}} \big){\re }^{\pm 2\ri\, \psi }=
q\hbar \sin \theta \cos 2\psi \pm \ri \,q\hbar \sin \theta \sin 2\psi =
{\widetilde{\pi }}_{1}\pm \ri \, {\widetilde{\pi }}_{2}
\end{equation*}%
do. Moreover, we have 
\begin{equation*}
\{{\widetilde{\pi }}_{\pm },{\widetilde{\pi }}_{3}\}=\{{\widetilde{\pi }}_{1},
{\widetilde{\pi }}_{3}\}\pm \ri \,\{{\widetilde{\pi }}_{2},{\widetilde{\pi }}_{3}\}=
2{\widetilde{\pi }}_{2}\pm 2\ri \, {\widetilde{\pi }}_{1}=
\pm 2\ri ({\widetilde{\pi }}_{1}\pm \ri \, {\widetilde{\pi }}_{2})=
\pm 2\ri \,{\widetilde{\pi }}_{\pm }.
\end{equation*}

Under reduction the quantum operator $\ttfrac{1}{\ri \hbar} \, {\mathbf{Q}}_{E}$ corresponds to the reduced
quantum operator $\ttfrac{1}{\ri \hbar}\, {\mathbf{Q}}_{{\widetilde{\pi }}_{4}}$ and the quantum operator 
$\ttfrac{1}{\ri \hbar} \, {\mathbf{Q}}_{L}$ corresponds 
to the reduced quantum operator $\ttfrac{1}{\ri \hbar} \, {\mathbf{Q}}_{{\widetilde{\pi }}_{3}}$. \medskip

As in the discussion of quantization of coadjoint in \cite{cushman-sniatycki12a}, we can define 
shifting operators ${\mathbf{Q}}_{{\widetilde{\pi }}_{\pm }}$ as follows. Set 
\begin{equation*}
{\mathbf{Q}}_{{\widetilde{\pi }}_{-}}{\widetilde{\mathbf{e}}}_{p,q}=b_{p}\,{\widetilde{\mathbf{e}}}_{p-2,q}\mspace{10mu}\mathrm{and}\mspace{10mu} 
{\mathbf{Q}}^{\dagger}_{{\widetilde{\pi }}_{-}}{\widetilde{\mathbf{e}}}_{p,q} 
= b_{p+2}\, {\widetilde{\mathbf{e}}}_{p+2,q}
\end{equation*}%
and 
\begin{equation*}
{\mathbf{Q}}_{{\widetilde{\pi }}_{+}}{\widetilde{\mathbf{e}}}_{p,q}= 
b_{p+2}\,{\widetilde{\mathbf{e}}}_{p+2,q} \mspace{10mu}\mathrm{and}\mspace{10mu} 
{\mathbf{Q}}^{\dagger}_{{\widetilde{\pi }}_{+}}{\widetilde{\mathbf{e}}}_{p,q} = 
b_p {\widetilde{\mathbf{e}}}_{p,q},
\end{equation*}%
where $b_{p}\in \R $ and $b_{-q}=0=b_{q+2}$. Then 
\begin{equation*}
{\mathbf{Q}}_{{\widetilde{\pi }}_{+}}{\mathbf{Q}}_{{\widetilde{\pi }}_{-}}{\widetilde{\mathbf{e}}}_{p,q}
=b_p{\mathbf{Q}}_{{\widetilde{\pi }}_{+}}{\widetilde{\mathbf{e}}}_{p+2,q}
=b_p^{2}\,{\widetilde{\mathbf{e}}}_{p,q}\mspace{10mu}\text{and}%
\mspace{10mu} {\mathbf{Q}}_{{\widetilde{\pi }}_{-}}{\mathbf{Q}}_{{\widetilde{\pi }}_{+}}{\widetilde{\mathbf{e}}}_{p,q}=b_{p+2}^{2}\,{\widetilde{\mathbf{e}}}_{p,q}.
\end{equation*}%
So 
\begin{equation*}
\lbrack {\mathbf{Q}}_{{\widetilde{\pi }}_{+}},{\mathbf{Q}}_{{\widetilde{\pi }%
}_{-}}]\, {\widetilde{\mathbf{e}}}_{p,q}={\mathbf{Q}}_{{\widetilde{\pi }}_{+}}{\mathbf{Q}%
}_{{\widetilde{\pi }}_{-}}\,{\widetilde{\mathbf{e}}}_{p,q}-{\mathbf{Q}}_{{\widetilde{%
\pi }}_{-}} {\mathbf{Q}}_{{\widetilde{\pi }}_{+}}\,{\widetilde{\mathbf{e}}}_{p,q}
=(b_{p}^{2}-b_{p+2}^{2}){\widetilde{\mathbf{e}}}_{p,q}.
\end{equation*}%
From 
\begin{equation*}
\{{\widetilde{\pi }}_{+},{\widetilde{\pi }}_{-}\}=\{{\widetilde{\pi }}_{1}+\ri \, {\widetilde{\pi }}_{2},{\widetilde{\pi }}_{1}- 
\ri \, {\widetilde{\pi }}_{2}\}=-2\ri \,\{{\widetilde{\pi }}_{1},{\widetilde{\pi }}_{2}\}
=-4\ri \,{\widetilde{\pi }}_{3}
\end{equation*}%
it follows that we should have 
\begin{equation}
\lbrack {\mathbf{Q}}_{{\widetilde{\pi }}_{+}},{\mathbf{Q}}_{{\widetilde{\pi }%
}_{-}}]=\ri \hbar \,{\mathbf{Q}}_{\{{\widetilde{\pi }}_{+},{\widetilde{\pi }}%
_{-}\}}=\ri \hbar \,{\mathbf{Q}}_{-4\ri \,{\widetilde{\pi }}_{3}}=
4\hbar \,{\mathbf{Q}}_{{\widetilde{\pi }}_{3}}. 
 \label{eq-s8one}
\end{equation}%
Evaluating both sides of equation (\ref{eq-s8one}) on ${\widetilde{\mathbf{e}}}_{p,q}$
gives 
\begin{equation*}
b_{p+2}^{2}-b_{p}^{2}=-4{\hbar }^{2}p,
\end{equation*}%
for every $p=-q+2k$, where $0\leq k\leq q $ and $k\in \Z $. In
particular we have $4{\hbar }^{2}q=b_{q+2}^{2}-b_{q}^{2}=-b_{q}^{2}$, since 
$b_{q+2}=0$ and $4{\hbar }^{2}(-q)=b_{-q+2}^{2}-b_{-q}^{2}=b_{-q+2}^{2}$,
since $b_{-q}=0$. Now for every $p=-q+2k$, where $k\in \{0,1,2,\ldots ,q\}$
we have 
\begin{align}
b_{p}^{2}& =(b_{p}^{2}-b_{p-2}^{2})+(b_{p-2}^{2}-b_{p-4}^{2})+\cdots
+(b_{-q+2}^{2}-b_{-q}^{2})+b_{-q}^{2}  \notag \\
& =-4{\hbar }^{2}(p-2)-4{\hbar }^{2}(p-4)+\cdots -4{\hbar }^{2}(-q),\quad 
\text{since}\mspace{5mu }b_{-q}=0  \notag \\
& =-4{\hbar }^{2}\sum_{\ell =1}^{k}(p-2\ell ),\quad \text{since}\mspace{5mu}%
p=-q+2k  \notag \\
& =-4{\hbar }^{2}(pk-(k+1)k)=4{\hbar }^{2}k(q+1-k)  \notag \\
& =-{\hbar }^{2}(p+q)(q-p+2).  
\label{eq-s8onestar}
\end{align}%
Therefore 
\begin{displaymath}
\left\{ \begin{array}{rl}
b_p & \hspace{-5pt} = 2\hbar \sqrt{k(q+1-k)} = \hbar \sqrt{(p+q)(q-p+2)} \\
b_{p+2} & \hspace{-3pt} = 2\hbar \sqrt{(k+1)(q-k)} = \hbar \sqrt{(p+q)(q-p)} .
\end{array} \right. 
\end{displaymath}
So when $p\in \{-q,\,-q+2,\,-q+4,\ldots ,q-2,\,q\}$ and $2k=p+q$ we
have 
\begin{align}
{\mathbf{Q}}_{{\widetilde{\pi }}_{-}}{\widetilde{\mathbf{e}}}_{p,q}& 
=2 \hbar \sqrt{k(q+1-k)}\,{\widetilde{\mathbf{e}}}_{p-2,q} \notag \\
{\mathbf{Q}}_{{\widetilde{\pi }}_{+}}{\widetilde{\mathbf{e}}}_{p,q}& 
=2 \hbar \sqrt{(k+1)(q-k)}\,{\widetilde{\mathbf{e}}}_{p+2,q} 
\notag 
\end{align}%
We now show that equation (\ref{eq-s8one}) holds. We compute 
\begin{align}
({\mathbf{Q}}_{{\widetilde{\pi }}_{+}}{\mathbf{Q}}_{{\widetilde{\pi }}_{-}}-{%
\mathbf{Q}}_{{\widetilde{\pi }}_{-}}{\mathbf{Q}}_{{\widetilde{\pi }}_{+}}){%
\widetilde{\mathbf{e}}}_{p,q}& =(b_{p}^{2}-b_{p+2}^{2}){\widetilde{\mathbf{e}}}_{p,q} \notag \\
&\hspace{-1.5in} =-4{\hbar }^{2} \big((k+1)(q-k)-k(q+1-k)\big){\widetilde{\mathbf{e}}}_{p,q}  \notag \\
& \hspace{-1.5in}=4{\hbar }^{2}(-q+2k){\widetilde{\mathbf{e}}}_{p,q}=
4{\hbar }^{2}p\,{\widetilde{\mathbf{e}}}_{p,q} 
=4{\hbar }\,{\mathbf{Q}}_{{\widetilde{\pi }}_{3}}{\widetilde{\mathbf{e}}}_{p,q}.  \notag
\end{align}

Since ${\widetilde{\pi }}_{1}=\onehalf ({\widetilde{\pi }}_{+}+{\widetilde{\pi }}_{-})$ and 
${\widetilde{\pi }}_{2}=\ttfrac{1}{2\ri } \, ({\widetilde{\pi }}_{+}-{\widetilde{\pi }}_{-})$, we
get 
\begin{equation}
{\mathbf{Q}}_{{\widetilde{\pi }}_{1}}{\widetilde{\mathbf{e}}}_{p,q}= %
\onehalf  \, {\mathbf{Q}}_{{\widetilde{\pi }}_{+}}{\widetilde{\mathbf{e}}}_{p,q}+ 
\onehalf \, {\mathbf{Q}}_{{\widetilde{\pi }}_{-}}{\widetilde{\mathbf{e}}}_{p,q} =%
\onehalf b_{p+2}\,{\widetilde{\mathbf{e}}}_{p+2,q} + 
\onehalf b_p \,{\widetilde{\mathbf{e}}}_{p-2,q}  \notag
\end{equation}%
and similarly 
\begin{displaymath}
{\mathbf{Q}}_{{\widetilde{\pi }}_{2}}{\widetilde{\mathbf{e}}}_{p,q}= %
\ttfrac{1}{2\ri }\,  b_{p+2} \,{\widetilde{\mathbf{e}}}_{p+2,q} 
- \ttfrac{1}{2\ri } \, b_{p}\,{\widetilde{\mathbf{e}}}_{p-2,q}.  
\end{displaymath}
The following calculation shows that the linear operator ${\mathbf{Q}}_{{\widetilde{\pi}}_1}$ is self adjoint on 
${\widetilde{\mathfrak{H}}}_q$. 
\begin{align}
{\mathbf{Q}}^{\dagger }_{{\widetilde{\pi }}_1}{\widetilde{\mathbf{e}}}_{p,q} & = {%
\mathbf{Q}}^{\dagger }_{{\widetilde{\pi }}_{+}}{\widetilde{\mathbf{e}}}_{p,q} + {%
\mathbf{Q}}^{\dagger }_{{\widetilde{\pi }}_{-}}{\widetilde{\mathbf{e}}}_{p,q} =
b_{p+2}\, {\widetilde{e}}_{p+2,q} + b_p \, {\widetilde{\mathbf{e}}}_{p-2,q} = {%
\mathbf{Q}}_{{\widetilde{\pi }}_1}{\widetilde{\mathbf{e}}}_{p,q}.  \notag
\end{align}
Similarly, the operators ${\mathbf{Q}}_{{\widetilde{\pi }}_k}$ for $k=2,3,4$
are self adjoint. \medskip 

The operators ${\mathbf{Q}}_{{\widetilde{\pi }}_{1}}$, ${\mathbf{Q}}_{{%
\widetilde{\pi }}_{2}}$, and ${\mathbf{Q}}_{{\widetilde{\pi }}_{3}}$ satisfy
the commutation relations 
\begin{equation*}
\mbox{{\footnotesize $\lbrack {\mathbf{Q}}_{{\widetilde{\pi}}_{1}},
{\mathbf{Q}}_{{\widetilde{\pi}}_{2}}]=\ri \hbar \,(2{\mathbf{Q}}_{{\widetilde{\pi}}_{3}}),\mspace{10mu}
\lbrack {\mathbf{Q}}_{{\widetilde{\pi}}_{2}},{\mathbf{Q}}_{{\widetilde{\pi}}_{3}}]=\ri \hbar
\,(2Q_{{\widetilde{\pi}}_{1}}),\mspace{10mu}\lbrack {\mathbf{Q}}_{{\widetilde{\pi}}_{3}},
{\mathbf{Q}}_{{\widetilde{\pi}}_{1}}]=\ri \hbar \,(2{\mathbf{Q}}_{{\widetilde{\pi}}_{2}}) $}} ,
\end{equation*}%
because 
\begin{align}
\lbrack {\mathbf{Q}}_{{\widetilde{\pi }}_{1}},{\mathbf{Q}}_{{\widetilde{\pi }%
}_{2}}]{\widetilde{\mathbf{e}}}_{p,q}& =
{\mathbf{Q}}_{{\widetilde{\pi }}_{1}}{\mathbf{Q}}_{{\widetilde{\pi }}_{2}}{\widetilde{\mathbf{e}}}_{p,q}
-{\mathbf{Q}}_{{\widetilde{\pi }}_{2}}{\mathbf{Q}}_{{\widetilde{\pi }}_{1}}{\widetilde{\mathbf{e} }}_{p,q}  \notag \\
& \hspace{-0.5in}= \ttfrac{1}{2\ri } b_{p+2}\,{\mathbf{Q}}_{{\widetilde{\pi }}_{1}}
{\widetilde{\mathbf{e}}}_{p+2,q}-%
\ttfrac{1}{2\ri }b_{p}\,{\mathbf{Q}}_{{\widetilde{\pi }}_{1}}
{\widetilde{\mathbf{e}}}_{p-2,q}  \notag \\
& - \ttfrac{1}{2\ri }b_{p+2}\, {\mathbf{Q}}_{{%
\widetilde{\pi }}_{2}}{\widetilde{\mathbf{e}}}_{p+2,q}-
\ttfrac{1}{2\ri }b_{p}\, {\mathbf{Q}}_{{\widetilde{\pi }}_{2}}
{\widetilde{\mathbf{e}}}_{p-2,q} \notag \\
& \hspace{-0.5in}= \ttfrac{1}{2\ri } b_{p+2}
(\onehalf b_{p+4}\, {\widetilde{\mathbf{e}}}_{p+4,q}
+\onehalf b_{p+2}\, {\widetilde{\mathbf{e}}}_{p,q})  \notag \\
& -\ttfrac{1}{2\ri } b_{p}
(\onehalf b_{p}\,{\widetilde{\mathbf{e}}}_{p,q} 
+ \onehalf b_{p-2}\,{\widetilde{\mathbf{e}}}_{p-4,q} )  \notag \\
& -\onehalf b_{p+2}
(\ttfrac{1}{2\ri }b_{p+4}\,{\widetilde{\mathbf{e}}}_{p+4,q}
-\ttfrac{1}{2\ri }b_{p+2}\,{\widetilde{\mathbf{e}}}_{p,q})  \notag \\
& -\onehalf b_{p}
(\ttfrac{1}{2\ri }b_{p}\,{\widetilde{\mathbf{\mathbf{e}}}}_{p,q}
-\ttfrac{1}{2\ri }b_{p-2}\, {\widetilde{\mathbf{e}}}_{p-4,q})  \notag \\
& \hspace{-0.5in}=-\ttfrac{1}{4\ri } (2b_{p}^{2}-2b_{p+2}^{2}){\widetilde{\mathbf{e}}}_{p,q}
= \ri \hbar \,(2p\hbar ){\widetilde{\mathbf{e}}}_{p,q}  
=\ri \hbar (2{\mathbf{Q}}_{{\widetilde{\pi }}_{3}}){\widetilde{\mathbf{e}}}_{p,q}. \notag
\end{align}

Similarly, 
\begin{equation*}
[{\mathbf{Q}}_{{\widetilde{\pi }}_{2}},{\mathbf{Q}}_{{\widetilde{\pi }}_{3}}]%
{\widetilde{\mathbf{e}}}_{p,q}= \ri \hbar \,(2{\mathbf{Q}}_{{\widetilde{\pi }}_{1}}){\widetilde{\mathbf{e}}}_{p,q} \mspace{10mu}\mathrm{and} \mspace{10mu} [{\mathbf{Q}}_{{\widetilde{\pi }}_{1}},
{\mathbf{Q}}_{{\widetilde{\pi }}_{3}}]{\widetilde{\mathbf{e}}}_{m,q}= 
\ri \hbar \,(-2{\mathbf{Q}}_{{\widetilde{\pi }}_{2}}){\widetilde{\mathbf{e}}}_{p,q}.
\end{equation*}
Therefore the skew hermitian operators $\ttfrac{1}{\ri \hbar} \, {\mathbf{Q}}_{{\widetilde{\pi }}_{1}}$, 
$\ttfrac{1}{\ri \hbar} \, {\mathbf{Q}}_{{\widetilde{\pi }}_{2}}$, and 
$\ttfrac{1}{\ri \hbar}\, {\mathbf{Q}}_{{\widetilde{\pi }}_{3}}$ satisfy the bracket
relations 
\begin{equation*}
\mbox{\footnotesize $[\ttfrac{1}{\ri \hbar}\, {\mathbf{Q}}_{{\widetilde{\pi }}_1}, 
\ttfrac{1}{\ri \hbar}\, {\mathbf{Q}}_{{\widetilde{\pi }}_2}] =
\ttfrac{1}{\ri \hbar} (2{\mathbf{Q}}_{{\widetilde{\pi }}_3}), \, \, 
[\ttfrac{1}{\ri \hbar}\, {\mathbf{Q}}_{{\widetilde{\pi }}_1}, 
\ttfrac{1}{\ri \hbar}\, {\mathbf{Q}}_{{\widetilde{\pi }}_3}] =
\ttfrac{1}{\ri \hbar}(-2{\mathbf{Q}}_{{\widetilde{\pi }}_2}), \, \, 
[\ttfrac{1}{\ri \hbar}\, {\mathbf{Q}}_{{\widetilde{\pi }}_2}, 
\ttfrac{1}{\ri \hbar}\, {\mathbf{Q}}_{{\widetilde{\pi }}_3}] =
\ttfrac{1}{\ri \hbar}(2{\mathbf{Q}}_{{\widetilde{\pi }}_1}) $},
\end{equation*}%
which defines the Lie algebra $\su  (2)$. Compare with table 1. 

Next we construct a representation of $\su (2)$ by skew hermitian linear operators. The bracket relations
show that for $j=1,2,3$ the mapping 
\begin{equation*}
{\widetilde{\mu}}_{q+1}:\su  (2) \rightarrow \mathfrak{u} ({\widetilde{\mathfrak{H}} }_{q}, \mathbb{C}): 
\ttfrac{1}{\ri \hbar}\, {\mathbf{Q}}_{{\widetilde{\pi}}_j} \rightarrow 
{\widetilde{\mu}}_{q+1}( %
\ttfrac{1}{\ri \hbar}\, {\mathbf{Q}}_{{\widetilde{\pi }}_j}),
\end{equation*}
where ${\widetilde{\mu }}_{q+1}(\ttfrac{1}{\ri \hbar} \, Q_{{\widetilde{\pi }}_{j}}):{\widetilde{\mathfrak{H}}}_{q}\rightarrow {\widetilde{\mathfrak{H}}}_{q}:{\widetilde{\mathbf{e}}}_{p,q}\mapsto 
\ttfrac{1}{\ri \hbar}\, Q_{{\widetilde{\pi }}_{j}}{\widetilde{\mathbf{e}}}_{p,q}$ is a faithful complex $2q+1$-dimensional
representation of $\su  (2)$ on ${\widetilde{\mathfrak{H}}}_{q}$ by skew hermitian linear operators. This representation can be decomposed into a sum of two irreducible representations as
follows. Let 
\begin{align}
{\widetilde{\mathfrak{H}}}_{q}^{0} &=\bigoplus\limits_{p=0}^{q}\mathrm{span} 
\{({\mathbf{Q}}_{{\widetilde{\pi }}_{-}})^{p}\,  {\widetilde{\mathbf{e}}}_{q,q}\}=
\bigoplus\limits_{p=0}^{q}\mathrm{span} \{ \widetilde{\mathbf{e}}_{q-2p,q}\}  \notag
\end{align}
and 
\begin{align}
{\widetilde{\mathfrak{H}}}_{q}^{1} &=\bigoplus\limits_{p=0}^{q-1}\mathrm{span} 
\{({\mathbf{Q}}_{{\widetilde{\pi }}_{-}})^{p}\, {\widetilde{\mathbf{e}}}_{q-1,q}\}=
\bigoplus\limits_{p=0}^{q-1}\mathrm{span} 
\{ \widetilde{\mathbf{e}}_{q-1-2p,q}\}.  \notag
\end{align}%
Then, ${\widetilde{\mathfrak{H}}}_{q}^{0}$ and ${\widetilde{\mathfrak{H}}}_{q}^{1}$ are 
irreducible representations of $\su (2) $ of dimension $q+1$ and $q$, respectively. Moreover, 
${\widetilde{\mathfrak{H}}}_{q}={\widetilde{\mathfrak{H}}}_{q}^{0}\oplus {\widetilde{\mathfrak{H}}}_{q}^{1}$. Thus, Bohr-Sommerfeld-Heisenberg quantization of a coadjoint orbit of $\mathrm{SU}(2)$ leads to a 
\emph{reducible} $\su (2)$ representation by skew hermitian operators.
\par \hspace{.5in}\begin{tabular}{l}
\vspace{.05in} \\
\includegraphics[width=200pt]{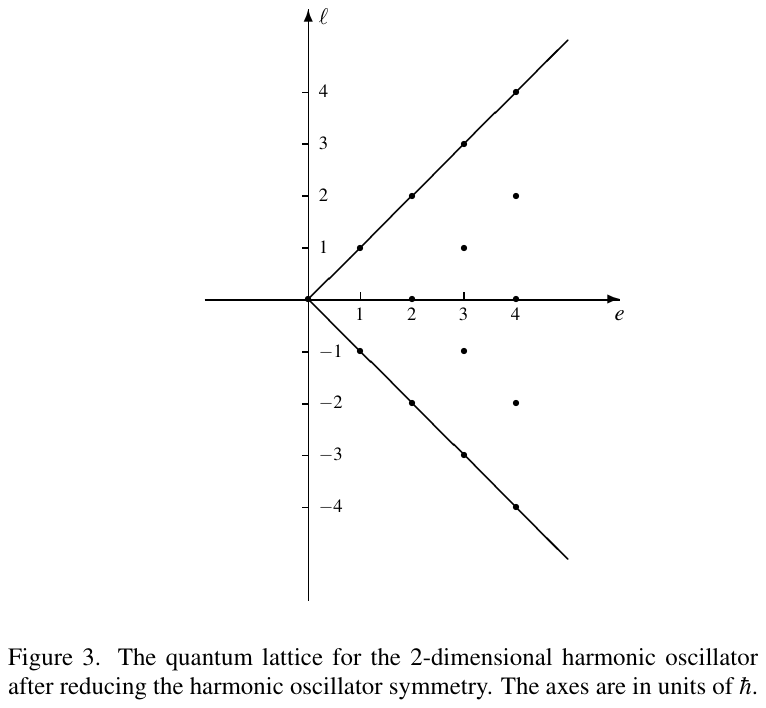} \\ 
\end{tabular}  

Using the Campbell--Baker--Hausdorff formula and the fact that $\mathrm{SU}(2)$ is compact, 
we can exponentiate the representation ${\widetilde{\mu }}_{q+1}$ to a unitary representation $R_{q}$
of $\mathrm{SU}(2)$ on ${\widetilde{\mathfrak{H}}}_{q}$. For more details see \cite{cushman-sniatycki12a}. 

\section{The fully reduced harmonic oscillator}

This section shows that the reason that, after reduction of the oscillator symmetry, the quantization of the 
Liouville integrable system on the $2$-sphere gives a 
\emph{reducible} representation of $\mathrm{SU}(2)$, is that the orginal $2$-degree of freedom harmonic oscillator has a ${\Z }_2$ symmetry. When this 
${\Z }_2$ symmetry is reduced, the quantized fully reduced system gives an \emph{irreducible} 
representation of $\mathrm{SU}(2)$.

\subsection[]{The classical fully reduced system}

In this subsection we reduce the \emph{full} symmetry group ${\Z }_2 \times S^1$ of the 
$2$-dimensional harmonic oscillator. \medskip 

Here the $S^1$ is generated by the flow of the 
harmonic oscillator vector field $X_{E}$ (\ref{eq-sec1dot}) on $T^{\ast }{\R }^4 = {\R }^4$ and 
the ${\Z }_2$-symmetry is generated by
\begin{equation}
\zeta :{\R }^4 \rightarrow {\R }^4:(x_1, x_2,y_1, y_2) \longmapsto (y_1, y_2,x_1,x_2) .
\label{eq-sec7onenw}
\end{equation}
This ${\Z }_2$-symmetry preserves the harmonic oscillator Hamiltonian E 
(\ref{eq-sec1dagger}), but is anti-symplectic as ${\zeta }^{\ast }\omega = \rd x_1 \wedge \rd y_1 + 
\rd x_2 \wedge \rd y_2 = - \omega $.  \medskip 

We now find the induced ${\Z }_2$ action on the generators ${\pi }_i$, $1\le i \le 4$ 
(\ref{eq-sec1invariants}) of the invariants of the oscillator symmetry. In 
terms of $(x,y)$ the $(\xi , \eta )$ coordinates (\ref{eq-sec2newone}) are 
\begin{displaymath}
(\xi ,\eta ) = \ttfrac{1}{\sqrt{2}}  \big( (x_1+y_2), (-x_1+y_2), (-x_2+y_1), (-x_2-y_1) \big) .
\end{displaymath} 
So the induced ${\Z }_2$ action on $(\xi , \eta )$ is the mapping 
\begin{displaymath}
\zeta : ( {\xi }_1, {\xi }_2, {\eta }_1, {\eta }_2) \mapsto (-{\eta }_2, -{\eta }_1, -{\xi }_2, -{\xi }_1 ).
\end{displaymath} 
Consequently, the induced action on the invariants ${\pi }_i$, $1\le i \le 4$ is the mapping 
\begin{displaymath}
\zeta : ({\pi }_1, {\pi }_2, {\pi }_3,{\pi }_4) \mapsto ({\pi }_1, {\pi }_2, -{\pi }_3, {\pi }_4).
\end{displaymath} 
Hence the algebra of ${\Z }_2 \times S^1$ invariant polynomials is generated by 
${\tau }_1 = {\pi }_1$, ${\tau }_2 = {\pi }_2$, ${\tau }_3 = {\pi }^2_3$, and ${\tau }_4 = {\pi }_4$ 
subject to the relation 
\begin{equation}
{\tau }_3 = {\pi }^2_3 = {\pi }^2_4 - {\pi }^2_1 -{\pi }^2_2 = {\tau }^2_4 - {\tau }^2_1-{\tau }^2_2, \, \, \, 
{\tau }_3 \ge 0 \, \, \& \, \, {\tau }_4 \ge 0,  
\label{eq-sec7defrel1}
\end{equation}
which defines the orbit space ${\R }^4/({\Z }_2 \times S^1)$ of the full symmetry 
group of the $2$-dimensional harmonic oscillator. \medskip 

On ${\R }^2\times {\R }^2_{\ge 0}$ 
we use the \emph{nonstandard} differential structure generated by the continuous functions 
${\tau }_1$, ${\tau}_2$, $\sqrt{{\tau }_3}$, and ${\tau }_4$. In other words, a function $f \in 
C^{\infty}({\R }^2\times {\R }^2_{\ge 0})$ if and only if there is a function 
$F \in C^{\infty}({\R }^4)$ such that $f = F({\tau }_1, {\tau}_2, \sqrt{{\tau }_3}, {\tau }_4)$, 
see \cite{cushman-bates}. Using table 1 in section 2, we see that the Poisson 
bracket ${\{ \, \, , \, \, \} }_{{\R }^2 \times {\R }^2_{\ge 0}}$ on 
$C^{\infty}({\R }^2\times {\R }^2_{\ge 0})$ has the structure matrix \medskip  

\noindent \hspace{1.15in}\begin{tabular}{c|cccc}
${\{{\tau }_i,{\tau  }_j\} }_{{\R }^2 \times {\R }^2_{\ge 0}}$ & ${\tau  }_1$ & ${\tau  }_2$ & $\sqrt{{\tau  }_3}$ & ${\tau  }_4$
\\ \hline
${\tau  }_1$ & $0$ & $2\sqrt{{\tau  }_3}$ & $-2{\tau  }_2$ & $0$ \\ 
${\tau  }_2$ & $-2{\tau }_2$ & $0$ & $2{\tau  }_1 $ & $0$ \\ 
$\sqrt{{\tau  }_3}$ & $2{\tau}_2 $ & $-4{\tau  }_1$ & $0$ & $0$ \\ 
${\tau }_4$ & $0$ & $0$ & $0$ & $0$
\end{tabular}
\vspace{.2in}

\noindent \hspace{.85in}\parbox[t]{3.5in}{Table 2. The structure matrix of the Poisson bracket 
${ \{ \, \, , \, \, \} }_{{\R }^2 \times {\R }^2_{\ge 0}}$ on 
$C^{\infty}({\R }^2 \times {\R }^2_{\ge 0})$.}\medskip 

\noindent The \emph{fully reduced} phase space ${\mathrm{H}}^{+}_e = S^2_e/{\Z }_2=   
(E^{-1}(e)/S^1)/{\Z }_2 \subseteq {\R }^3$ with coordinates $( {\tau }_1, {\tau }_2, {\tau }_3)$, which 
is also $E^{-1}(e)/$ $({\Z }_2 \times S^1)$, is defined by (\ref{eq-sec7defrel1}) 
and ${\tau }_4 =e$, that is, ${\tau }_3 = e^2 - {\tau }^2_1 - {\tau }^2_2$ with ${\tau }_3 \ge 0$. 
Geometrically, ${\mathrm{H}}^{+}_e$ is the upper hemisphere of the $2$-sphere $S^2_e$. Because 
${\tau }_4$ is a Casimir of the Poisson algebra $(C^{\infty}({\mathbf{R}}^2 \times {\mathbf{R}}^2_{\ge 0}), 
{\{ \, \, , \, \, \} }_{{\R }^2 \times {\R }^2_{\ge 0}}, \cdot )$, the Poisson 
bracket ${\{ \, \, , \, \, \} }_{H^{+}_e}$ on the space $C^{\infty}(H^{+}_e)$ of smooth functions on 
the fully reduced phase space $H^{+}_e$ has structure matrix given in table 3. Note that 
the structure matrix of the \linebreak 

\par\noindent \hspace{1.35in}\begin{tabular}{c|ccc}
${\{{\tau }_i,{\tau  }_j\} }_{H^{+}_e}$ & ${\tau  }_1$ & ${\tau  }_2$ & $\sqrt{{\tau  }_3}$ 
\\ \hline
\rule{8pt}{0pt}${\tau  }_1$ & $0$ & $2\sqrt{{\tau  }_3}$ & $-2{\tau  }_2$  \\ 
${\tau  }_2$ & $-2 \sqrt{{\tau}_3}$ & $0$ & $2{\tau  }_1 $  \\ 
$\sqrt{{\tau  }_3}$ & $2{\tau }_2 $ & $-2{\tau  }_1$ & $0$  \\ 
\end{tabular}
\vspace{.2in}
\par \noindent \hspace{.85in}\parbox[t]{3.5in}{Table 3. The structure matrix of the Poisson bracket 
${\{ \, \, , \, \, \} }_{H^{+}_e}$ on $C^{\infty}(H^{+}_e)$.}\medskip

\noindent Poisson bracket ${\{ \, \, , \, \, \} }_{H^{+}_e}$ in table 3 defines a 
Lie algebra, which is isomorphic to $\su (2)$. The differential structure $C^{\infty}(H^{+}_e)$ is 
generated by the functions ${\tau }_1$, ${\tau }_2$, and $\sqrt{{\tau }_3}$. \medskip 

The next argument shows that the differential space $(H^{+}_e, C^{\infty}(H^{+}_e))$ 
is locally diffeomorphic to a closed subset of Euclidean space. In other words, the differential space 
$(H^{+}_e, C^{\infty}(H^{+}_e))$ is subcartesian, see \cite{sniatycki12}. The image under the restriction of the injective map 
\begin{displaymath}
\phi :{\R }^2 \times {\R }_{\ge 0} \rightarrow {\R }^2 \times {\R }_{\ge 0}:
({\tau }_1, {\tau}_2, \sqrt{{\tau }_3}) \longmapsto ({\pi }_1, {\pi }_2, {\pi }_3) = 
\big( {\tau }_1, {\tau }_2, (\sqrt{{\tau }_3})^2 \big) 
\end{displaymath}
to $H^{+}_e$ is $S^2_e \cap \{ {\pi }_3 \ge 0 \} $. So  
${\phi }^{\ast }\big( C^{\infty}(S^2_e \cap \{ {\pi }_3 \ge 0 \}) \big) = 
C^{\infty}(H^{+}_e)$. Thus, the differential spaces $(H^{+}_e, C^{\infty}(H^{+}_e))$ and 
$(S^2_e \cap \{ {\pi }_3 \ge 0 \} , C^{\infty}(S^2_e \cap \{ {\pi }_3 \ge 0 \}))$ are diffeomorphic. 
Because $S^2_e \cap \{ {\pi }_3 \ge 0\} $ is a semialgebraic variety in ${\R }^3$, being 
defined by ${\pi }^2_1 +{\pi }^2_2+{\pi }^2_3 = e^2$ and ${\pi }_3 \ge 0$, it is a closed subset of 
${\R }^3$. Hence $C^{\infty}(S^2_e \cap \{ {\pi }_3 \ge 0 \})  = 
C^{\infty}({\R }^3)|(S^2_e \cap \{ {\pi }_3 \ge 0 \} )$. This implies that $(S^2_e \cap \{ {\pi }_3 \ge 0 \} , C^{\infty}(S^2_e \cap \{ {\pi }_3 \ge 0 \}))$ is a locally compact subcartesian differential space, see \cite[Chapter VII]{cushman-bates}.  \medskip 

The fully reduced Hamiltonian of the fully reduced system is 
\begin{equation}
L^{\scriptscriptstyle \vee }_e: H^{+}_e \subseteq {\R }^2 \times {\R }_{\ge 0} \rightarrow 
\R : ({\tau }_1, {\tau }_2, {\tau }_3) \longmapsto \sqrt{{\tau }_3} 
\label{eq-sec7ham}
\end{equation}
and the fully reduced equations of motion on ${\R }^3$ are 
\begin{align}
{\dot{\tau }}_1 & = {\{ {\tau }_1, L^{\scriptscriptstyle \vee}_e \}}_{H^{+}_e} 
= {\{ {\tau }_1, {\tau  }_3 \} }_{H^{+}_e} = -2{\tau }_2 \notag \\
{\dot{\tau }}_2 & = {\{ {\tau }_2, L^{\scriptscriptstyle \vee }_e \}}_{H^{+}_e} 
= { \{ {\tau }_2, {\tau  }_2 \}}_{H^{+}_e} = 2{\tau }_1
\label{eq-sec7eqnsmot} \\
{\dot{\tau }}_3 & = { \{ {\tau }_3, L^{\scriptscriptstyle \vee }_e \} }_{H^{+}_e} =0. \notag 
\end{align}
Since $(\sqrt{{\tau }_3})^2+{\tau }^2_1+{\tau }^2_2$ is a Casimir of the Poisson algebra 
$(C^{\infty}({\R }^2 \times {\R }^2_{\ge 0}),$ \linebreak 
${\{ \, \, , \, \, \} }_{{\R }^2 \times {\R }^2_{\ge 0}}, \cdot )$, 
the fully reduced phase space $H^{+}_e$ is invariant under the flow of 
(\ref{eq-sec7eqnsmot}). So (\ref{eq-sec7eqnsmot}) restricted to $H^{+}_e$ gives the 
integral curves of the fully reduced Hamiltonian vector field $X_{L^{\scriptscriptstyle \vee }_e}$ on the 
differential space $\big( H^{+}_e, C^{\infty}(H^{+}_e) \big) $. \medskip

The map 
\begin{displaymath}
\begin{array}{l}
\psi : {\overline{D}}^2_e = \{  ( {\tau }_1, {\tau }_2 ) \in {\R }^2 \, | \, 
{\tau }^2_1 +{\tau }^2_2 \le e^2 \} \rightarrow 
H^{+}_e \subseteq {\R }^3: \\
\rule{0pt}{13pt} \hspace{.75in} ( {\tau }_1, {\tau }_2) \longmapsto ({\tau }_1, {\tau }_2, 
\sqrt{e^2-{\tau }^2_1 -{\tau }^2_2}) 
\end{array}
\end{displaymath}
is a diffeomorphism of the differential space $({\overline{D}}^2_e, C^{\infty}({\overline{D}}^2_e))$ 
onto the differential space $(H^{+}_e, C^{\infty}(H^{+}_e))$ provided we say that 
$f \in C^{\infty}({\overline{D}}^2_e)$ 
if and only if there is an $F \in C^{\infty}(H^{+}_e)$ such that $f = {\psi }^{\ast }F$. In other words, \linebreak 

\noindent \hspace{1in}\begin{tabular}{c|cc}
${\{{\tau }_i,{\tau  }_j \}} _{{\overline{D}}^2_e}$ & ${\tau  }_1$ & ${\tau  }_2$  
\\ \hline
\rule{0pt}{12pt}${\tau  }_1$ & $0$ & $2\sqrt{e^2-{\tau }^2_1-{\tau }^2_2}$   \\ 
${\tau  }_2$ & $-2\sqrt{e^2-{\tau }^2_1-{\tau }^2_2}$ & $0$   
\end{tabular}
\vspace{.2in}

\noindent \hspace{.8in}\parbox[t]{3.5in}{Table 4. The structure matrix of the Poisson bracket 
${\{ \, \, , \, \, \} }_{{\overline{D}}^2_e}$ on $C^{\infty}({\overline{D}}^2_e)$.}\medskip

\noindent $C^{\infty}({\overline{D}}^2_e) = {\psi }^{\ast }(C^{\infty}(H^{+}_e))$. Note that topologically 
${\overline{D}}^2_e$ is a closed $2$-disk, but its differential structure $C^{\infty}({\overline{D}}^2_e)$ is 
\emph{nonstandard}. The structure matrix of the Poisson bracket 
$\{ \, \, , \, \, \}|_{{\overline{D}}^2_e}$ on $C^{\infty}({\overline{D}}^2_e)$ is given in table 4. \medskip 

\noindent The fully reduced Hamiltonian on ${\overline{D}}^2_e$ is 
\begin{displaymath}
{\widehat{L}}_e = {\psi }^{\ast }L^{\scriptscriptstyle \vee}_e: {\overline{D}}^2_e \subseteq {\R }^2 
\rightarrow \R : ({\tau }_1, {\tau }_2) \longmapsto \sqrt{e^2 -{\tau }^2_1 -{\tau }^2_2}. 
\end{displaymath}
The Hamiltonian vector field $X_{{\widehat{L}}_e}$ of ${\widehat{L}}_e$ on the 
differential space $({\overline{D}}^2_e, C^{\infty}({\overline{D}}^2_e))$ has integral curves governed by 
\begin{align}
{\dot{\tau }}_1 & = {\{ {\tau }_1, {\widehat{L}}_e \} }_{{\overline{D}}^2_e} =
\mbox{$\frac{\scriptstyle 1}{\scriptstyle 2\sqrt{e^2 -{\tau }^2_1 -{\tau }^2_2}}$} \big( -2{\tau }_2 
{\{ {\tau }_1, {\tau }_2 \} }_{{\overline{D}}^2_e} \big)  = - 2{\tau }_2 \notag \\
{\dot{\tau }}_2 & = {\{ {\tau }_2, {\widehat{L}}_e \} }_{{\overline{D}}^2_e} = 2 {\tau }_1 . \notag 
\end{align}
 
The function $\widehat{A} = \onehalf \, \sqrt{e^2 -{\tau }^2_1 -{\tau }^2_2}$ is 
an action function on ${D}^{\ast }_e = {D}^2_e \setminus \{ (0,0) \}$, because its Hamiltonian vector field $X_{\widehat{A}}$ is $-{\tau }_2\frac{\partial }{\partial {\tau }_1}
+ {\tau }_1\frac{\partial }{\partial {\tau }_2}$ has integral curves $t \mapsto 
${\tiny $\begin{pmatrix} \cos t & -\sin t \\
\sin t & \cos t \end{pmatrix} \, \begin{pmatrix} {\tau }^0_1 \\ \rule{0pt}{8pt}{\tau }^0_2 \end{pmatrix} $} with 
$({\tau }^0_1, {\tau}^0_2) \in {D}^{\ast }_e$, which are periodic of period $2\pi $. The corresponding angle is 
$\widehat{\vartheta } =  {\tan }^{-1}\frac{{\tau }_2}{{\tau }_1}$. Using table 4 we see that the fully reduced symplectic form ${\widehat{\omega }}_{{D}^{\ast }_e}$ on 
${D}^{\ast }_e$ is $\frac{1}{2\sqrt{e^2 -{\tau }^2_1 -{\tau }^2_2}} \, \rd {\tau }_2 \wedge \rd {\tau }_1$, which equals 
$\rd \widehat{A} \wedge \rd \widehat{\vartheta }$.

\subsection[]{The quantized fully reduced harmonic oscillator}

In this subsection we give the Bohr-Sommerfeld-Heisenberg quantization of the fully reduced 
harmonic oscillator $({\widehat{L}}_e, {\overline{D}}^2_e, {\{ \, \, , \, \, \} }_{{\overline{D}}^2_e})$, where 
$e$ is fixed and is nonnegative.  \medskip 

A straightforward calculation shows that the volume ${\mathrm{vol}}_{D^{\ast }_e}$ of the punctured 
$2$-disk $D^{\ast }_e$ using the fully reduced symplectic form ${\widehat{\omega }}_{{D}^{\ast }_e}$ on 
${D}^{\ast }_e$ is $\frac{1}{2\sqrt{e^2 -{\tau }^2_1 -{\tau }^2_2}} \, \rd {\tau }_2 \wedge \rd {\tau }_1$ is 
$2\pi \, e$. Thus $({D}^{\ast }_e, {\widehat{\omega }}_{{D}^{\ast }_e})$ is quantizable if there is a positive 
integer $q$ such that $2\pi \, e = h \, q$, that is, $e = \hbar \, q$. Note that $q$ is a fixed positive 
integer. \medskip 

The $1$-torus ${\widehat{A}}^{-1}(\ell )$ on $D^{\ast }_e$ is a Bohr-Sommerfeld 
$1$-torus if there is an integer $m$ such that 
\begin{displaymath}
m \, h = \int_{{\widehat{A}}^{-1}(\ell )} \widehat{A}\, \rd \widehat{\vartheta } = \ell (2 \pi  ), 
\end{displaymath}
that is, $\ell = m \, \hbar $. The second equality above follows because ${\widehat{A}}^{-1}(\ell )$ 
is parametrized by $\widehat{\vartheta} \mapsto \sqrt{e^2-{\ell }^2}(\cos \widehat{\vartheta }, \sin \widehat{\vartheta })$ for $\widehat{\vartheta } \in [0, 2\pi ]$. On 
$D^{\ast }_e$ we have ${\ell }^2 =  e^2 - {\tau }^2_1 -{\tau }^2_2$, which implies $e^2 \ge {\ell }^2$, that is, 
$0 < m < q$. Thus there are $q+1$ Bohr-Sommerfeld tori on ${\overline{D}}^2_e$, namely,   
$q-1$ regular Bohr-Sommerfeld $1$-tori $\{ ({\tau }_1, {\tau }_2) \in D^{\ast }_e \, | \, 
{\tau }^2_1 +{\tau }^2_2 = (q^2- m^2) {\hbar }^2 \} $ with 
$0 < m < q$, one singular $1$-torus $\partial {\overline{D}}^2_e = 
\{ ({\tau }_1, {\tau }_2) \in D^{\ast }_e \, | \,  
{\tau }^2_1 +{\tau }^2_2 =e^2 \}$, and a singular point $\{ (0,0) \} $. \medskip 

Corresponding to each of the $q+1$ Bohr-Sommerfeld tori constructed in the preceding paragraph, we 
associate a vector ${\widehat{\mathbf{e}}}_m$, $0 \le m \le q$, which are orthogonal under the hermitian 
inner product $\langle {\widehat{\mathbf{e}}}_k, {\widehat{\mathbf{e}}}_{k'} \rangle = {\delta }_{k,k'}$. Their  
span is the Hilbert space ${\widehat{\mathfrak{H}}}_q$, which is the eigenspace of the self adjoint 
quantum operator ${\mathbf{Q}}_{\widehat{A}}$ corresponding to the eigenvalue $\hbar \, q$. 
The space of quantum states of the fully reduced harmonic oscillator is the Hilbert space $\widehat{\mathfrak{H}} = {\bigoplus}_{q \in {\Z }_{\ge 0}} {\widehat{\mathfrak{H}}}_q$. The shifting operators on 
${\widehat{\mathfrak{H}}}_q$ in the Bohr-Sommerfeld-Heisenberg quantization of the fully reduced harmonic oscillator are the lowering operator ${\mathbf{a}}_1({\widehat{\mathbf{e}}}_m) = 
\left\{ \right. ${\tiny $\begin{array}{cl} 
{\widehat{\mathbf{e}}}_{m-1}, & \mbox{if $1 \le m \le q$} \\
0, & \mbox{if $m=0$} \end{array}$} and the raising operator ${\mathbf{a}}_2({\widehat{\mathbf{e}}}_m) 
= {\widehat{\mathbf{e}}}_{m+1}$, if $0 \le m \le q-1$. These operators correspond to the quantum operators 
${\mathbf{Q}}_{{\re}^{-\ri \widehat{\vartheta }}}$ and ${\mathbf{Q}}_{{\re}^{\ri \widehat{\vartheta }}}$, respectively.  
Since ${\re}^{\pm \ri \widehat{\vartheta }}$ are not smooth functions on ${\overline{D}}^2_e$, we \linebreak 
replace them 
with the smooth functions ${\tau }_{\pm } = \ttfrac{1}{\ri \hbar} ( {\tau }_1 \pm \ri \, {\tau }_2) = 
\ttfrac{1}{\ri \hbar}\, r {\re}^{\pm \ri \widehat{\vartheta}}$, where 
$r = \sqrt{{\tau }^2_1 + {\tau }^2_2}$. Using table 4 we get   
\begin{displaymath}
{\{ {\tau }_{+}, {\tau }_{-} \} }_{{\overline{D}}^2_e} = -\ri \, {\{ {\tau }_1, {\tau }_2 \} }_{{\overline{D}}^2_e} 
= -2\ri \, \sqrt{e^2-{\tau }^2_1-{\tau }^2_2} = -4\ri \, \widehat{A}, 
\end{displaymath} 
${\{{\tau }_{+},  \widehat{A} \} }_{{\overline{D}}^2_e} =  \ri \,  {\tau }_{+}$, and  
${\{ {\tau }_{-}, \widehat{A} \} }_{{\overline{D}}^2_e} = - \ri \,  {\tau }_{-}$. Thus $\{ {\tau}_1, {\tau }_2, 
\widehat{A} \} $ span a Lie algebra under ${\{ \, \, , \, \,  \} }_{{\overline{D}}^2_e}$, which is isomorphic to 
$\su (2)$. Define quantum shifting operators on ${\widehat{\mathfrak{H}}}_q$ by 
${\mathbf{Q}}_{{\tau }_{-}}{\widehat{\mathbf{e}}}_m = c_m \, 
{\widehat{\mathbf{e}}}_{m-1}$ if $1 \le m \le q$ and $0$ if $m=0$ and ${\mathbf{Q}}_{{\tau }_{+}}{\widehat{\mathbf{e}}}_m= c_{m+1} \, {\widehat{\mathbf{e}}}_{m+1}$, if $0 \le m \le q-1$ and their adjoints 
${\mathbf{Q}}^{\dagger}_{{\tau }_{-}}{\widehat{\mathbf{e}}}_m= c_{m+1}\, {\widehat{\mathbf{e}}}_{m+1}$ and 
${\mathbf{Q}}^{\dagger}_{{\tau }_{+}}{\widehat{\mathbf{e}}}_m= c_{m-1}\, {\widehat{\mathbf{e}}}_{m-1}$ if 
$1\le m \le q$ and $0$ if $m=0$. Here $c_m \in \mathbb{C}$ for each $0 \le m \le q$. An 
argument analogous to the one used in section 5 determines the $q+1$ coefficients 
$c_m$ for $0 \le m \le q$. We get 
$c_m =${\tiny $\left\{ \begin{array}{cl} \ri \hbar \, \sqrt{2m(m-1)}, & \mbox{if $1 \le m \le q$} \\
0, & \mbox{if $m=0$.} \end{array} \right. $} Thus for each $q \in {\Z }_{\ge 0}$ the operators 
$\{ {\mathbf{Q}}_{{\tau }_1}, {\mathbf{Q}}_{{\tau }_2}, 
{\mathbf{Q}}_{\widehat{A}} \} $ give rise to a $q+1$-dimensional irreducible representation of 
$\su (2)$ on ${\widehat{\mathfrak{H}}}_q$. 

\section{Summary}
\label{sec7}

We summarize the contents of this paper. In sections 2---6 we give the 
Bohr-Sommerfeld-Heisenberg quantization of the following Liouville integrable systems: 
the two degree of freedom harmonic oscillator $(E, L, T^{\ast }{\R }^2, \omega )$; 
its reduction by the $S^1$ oscillator symmetry on a fixed energy level $e$ to the one 
degree of freedom system $({\widetilde{L}}_e, S^2_e, {\widetilde{\omega }}_e)$; and 
the reduction of the latter system under the ${\Z }_2$ symmetry to a one 
degree of freedom system $({\widehat{L}}_e, {\overline{D}}^2_e, {\{ \, \, , \, \, \} }_{{\overline{D}}^2_e} )$. 
In each of the above integrable systems we find the Hilbert space of quantum states 
using Bohr-Sommerfeld quantization. Using the lowering and raising operators coming from the Heisenberg quantization, we obtain skew symmetric quantum operators on three spaces of quantum states. Each of the 
three collections of quantum operators gives rise to a representation of a finite dimensional Lie algebras, namely,
${\mathcal{L}}_2$, $\su (2)$, and $\su (2)$, respectively. Each of these representations decompose their respective space of quantum states into: a countably infinite number of copies of an infinite dimensional irreducible representation of ${\mathcal{L}}_2$; two finite dimensional irreducible representations of $\su (2)$ of dimension $q$ and $q+1$; and a single finite dimensional irreducible representation of $\su (2)$ of dimension $q+1$.

\end{document}